\documentclass[12pt,a4paper]{article}
\usepackage{amsfonts,amsmath,amssymb}
\usepackage{graphicx,color}
\usepackage[export]{adjustbox}
\usepackage{rotating}
\usepackage{parskip}
\usepackage{hyperref}
\usepackage{tabularx}
\usepackage{enumitem}

\usepackage{chngcntr}

\counterwithin*{equation}{section}

\newcommand{\K}[3]{\left#1#3\right#2}
\newcommand{\s}{\hspace{2mm}}

\begin{document}

	\begin{center}
		{\LARGE\bf Multi-Dimensional Martingales from Mutual Information}
		\\[1cm]
		\large{Michael M. Kay}\footnote{e-mail: michael.kay@outlook.de}
		\\[2cm]
		
		\parbox{14cm}
		{ \small
		{\bf Abstract.} In the context of Risk Neutral Pricing theory, we consider the classic problem of calibrating a martingale over $\mathbb{R}^n$ to a finite number of marginals thereof, or more practically, to prices of an arbitrary finite set of (joint) European contingent claims. For $n=1$, one can rely on the work of Dupire, while for $n\geq 2$ an analogous natural unique construction seems to be lacking. We provide such a unique candidate as the result of pure Martingale Entropic Optimal Transport. As a byproduct, the latter allows us to obtain a constructive proof of a classic result of Strassen. Finally, and in contrast to the proposed approach, we prove a result that demonstrates how a certain class of local correlation models fails in general to calibrate to basket option prices, particularly in the foreign exchange market.
		}
	\end{center}
\tableofcontents

\setcounter{section}{-1}
\section{Introduction}

The aim of the present work is to introduce a concrete method to calibrate a Risk Neutral Markov process to an arbitrary set of European contingent claims\footnote{ These are trades whose payoffs have a single cashflow at maturity that depends solely on the value of the tradables at maturity.} that depend on multiple tradables. The archetypal example, where such a type of problem arises naturally is in the foreign exchange market. In particular, when pricing a trade that depends in some non-trivial manner, say on three currencies, which we shall refer to as $d$ (domestic), $f_1$, $f_2$, a minimal set of vanilla products one needs to calibrate to are European options on the FX rates $X_{d \leftarrow f_1}, X_{d \leftarrow f_2}, X_{f_1 \leftarrow f_2}$. Assuming no arbitrage, the latter is the ratio of the first two, hence its corresponding options can be recast as particular basket options on the first two. A popular class of models in this context is the one studied in \cite{Guyon13}, however such a class fails in general (see Appendix A) to calibrate to arbitrage free prices. Other popular approaches for such a type of problems involve ideas from martingale Optimal Transport (see the comprehensive review \cite{LabordereBook}, and e.g., \cite{LabordereMEOT}, \cite{GuyonVIX}, \cite{Guo1}, \cite{Guo2}). Such an approach is more general and has a more general scope than what we propose, in particular the products it allows to calibrate to are in principle arbitrary. However a specialization thereof that could potentially be viewed as an extension of the work by Dupire (\cite{Dupire}) in the case of a single tradable, seems to be lacking. When dealing with a single tradable, the work of Dupire, aside from regularity considerations, allows one to construct a unique continuous local martingale given, say, the collection of all European call option prices $C(t,K)$ for all maturities $t$ and strikes $K$. In the setting of zero interest rates, in particular, the local volatility of such a process is obtained as
\begin{align}
\sigma^2_t(K) = 2\frac{\partial_t C(t, K)}{\partial_K^2 C(t, K)}.
\end{align}

In practice, one cleans (for no-arbitrage) and interpolates/extrapolates the function $C(t,K)$ from a finite set of maturities and strikes. The above is nothing else than the result of solving for $\sigma_t^2(x)$ in the forward Kolmogorov (Fokker-Planck) equation
\begin{align}
\partial_t \rho_t(x) = \frac{1}{2}\partial_x^2 \K{(}{)}{\sigma^2_t(x)\rho_t(x)}
\end{align}

given the marginal probability densities $\rho_t$ over $\mathbb{R}$. The appeal of Dupire's result is that it offers a quick analytic solution (which however becomes more involved in the case of stochastic rates) and that when $\sigma^2_t$ exists, it is unique. The existence of $\sigma_t^2$ is essentially guaranteed in the absence of calendar arbitrage. More precisely, in the case of discrete time, when $\rho_s \preceq \rho_t$ for $s < t$ where $\preceq$ denotes \textsl{convex order}, a classic result by Strassen (see \cite{Strassen}) guarantees the existence of a Markov martingale with those marginals. Such a result is true not only in dimension one (over $\mathbb{R}$), but for Borel measures over a separable Banach space. In the case of continuous time the analogous result is due to Kellerer (see \cite{Kellerer}) but only in dimension one, while a generalization to $\mathbb{R}^n$, requiring smoothing, was obtained very recently in \cite{Schachermayer}. Both in discrete and continuous time, however, the results are not constructive enough to make them useful in practice. Moreover, such Markov processes are no longer unique in $\mathbb{R}^n$ for $n > 1$.  The method we propose stems from a constructive proof of Strassen's theorem for Radon probability measures over $\mathbb{R}^n$, whereby one obtains the unique Markov martingale that minimizes the mutual information of the joint distribution of the process\footnote{By mutual information we mean the Kullback-Leibler divergence of the joint distribution relative to the product of its marginals.}. One way of interpreting such a martingale, is as the most unbiased martingale given marginals free of calendar arbitrage. We will  refer to such a martingale as a \textsl{Hessian martingale}. Subsequently, we construct the unique Markov (Hessian) martingale that minimizes mutual information, where the constraints are not the knowledge of the full marginals, but rather the prices of a finite number of European Contingent claims. This latter result can be used as an effective calibration method in practice. The method we propose has the following properties in common with Dupire's:
\begin{itemize}
\item Up to numerical error, it succeeds if and only if there is no market arbitrage and can be also employed to detect/``clean away" arbitrage.
\item It is ``model independent" in the sense that the form of the resulting Markov process is completely and solely determined by the form of the payoffs. 
\end{itemize}

In contrast to Dupire's approach:
\begin{itemize}
\item The form of the volatility matrix (or covariance matrix) is not determined from a simple algebraic formula, but requires instead the solution of a convex optimization problem.
\item The proposed approach only needs an arbitrary finite number of European contingent claims. In particular it does not require ad-hoc interpolation/extrapolation of prices or implied volatility (hyper-)surfaces.
\end{itemize}

 It is essential to mention that the idea of looking for a (Markov) process that minimizes the relative entropy between the joint distribution of the process and a given, reference one ( in our case the product of its marginals ) is nothing new, and in the physics literature goes back at least to Schrödinger (see \cite{Schroedinger}), while in the mathematical finance literature such a line of research was pioneered in \cite{Avellaneda1, Avellaneda2, Avellaneda3}. However, one added value we claim to provide, relative to that line of research, lies in unravelling the full structure of the resulting process when specializing to European contingent claims. Moreover, contrary to such approaches, here the martingale property is imposed a priori.

The present note is subdivided into 7 sections plus 4 appendices and it has 3 main results called Theorem 1,2,3, the first of which is a version of the classic result by Strassen specialized to $\mathbb{R}^n$, for which, however, we offer a constructive proof\footnote{We note that for $n=1$ a novel proof in the same vein was also obtained recently in \cite{Wiesel}.}:
\begin{enumerate}
\setcounter{enumi}{0}
\item Revisitation of Strassen's theorem: Theorem 1.
\item Intuitive explanation of the construction of the unique ``most unbiased" Markov martingale incarnating Strassen's result.
\item Proof of Theorem 1. 
\item Intuitive understanding of the type of martingale Markov process constructed given the marginals.
\item Sketch of what the process looks like in continuous time. It shows that under mild assumptions, the process is a continuous local martingale, where the covariance matrix is the inverse of the Hessian of a time-dependent convex function on the space of tradables.
\item Statement and Proofs of Theorem 2,3: these generalize Theorem 1, by not requiring knowledge of the full marginals, but only of the prices of some European contingent claims. Theorem 2: shows exactly under what conditions the calibration procedure can be sequential. This condition incidentally is generically only satisfied in dimension 1. 

Theorem 3 shows how to construct the Hessian martingale to prices of an arbitrary finite set of European contingent claims and explains the full calibration approach. Pleasant features thereof are that there are exactly as many unknowns as there are calibration targets and the optimization is convex.
\item A sketch of the calibration methodology in the Monte Carlo framework.
\end{enumerate}

\begin{itemize}
\item[A.] We show how the class of local correlation models studied in \cite{Guyon13} fails in general to calibrate to arbitrage free basket option prices in dimension $2$ and higher. Specifically, our result shows that such a class is insufficient in general to calibrate simultaneously to arbitrage free prices of FX options relative to three currencies.
\item[B.] We present details of proofs of propositions involving convex functions, which are necessary for Theorem 2,3.
\item [C.] We show how Hessian martingales are optimal in the sense of mutual information, among all (not necessarily Markov) martingales. The latter in a sense explains the naturality of Markov processes when only the marginals of the process are known.
\item [D.] Simple examples revisited in the light of the formalism here explained. Useful to become familiar with the formalism. 
\end{itemize}

\textbf{Remark on notation:} We caution the reader on the potentially non-standard use of notation throughout the paper. In particular:
\begin{itemize}
\item Coordinate functions on $\mathbb{R}^n$, when indexed, have indices upstairs: $x^i$ as is customary in e.g. differential geometry / physics literature, but non-standard elsewhere.
\item Partial derivatives w.r.t. $x^i$ are denoted as $\partial_i$.
\item An expression of the form $\partial_i x^i$ implies summation over the index and thus would yield $n$ in this example for $x \in \mathbb{R}^n$.
\item An expression of the form $(1,x)$, where $x \in \mathbb{R}^n$, stands for an element of $\mathbb{R}^{n+1}$ whose first component is $1$ and the remaining $n$ components coincide with those of $x$.
\item A potentially uncomfortable, yet convenient notation that will be employed at times, is the expression $d\mu(x) = \rho(x)d^n x$ for a Radon probability measure $\mu$ over $\mathbb{R}^n$, where $d^n x$ refers to the Lebesgue measure, and where $\rho$ can be regarded as a ``generalized" density, in particular it will not necessarily be in $L^1(\mathbb{R}^n)$.
\item Covariance matrices are viewed as inverses $g^{-1}$ of Riemannian metrics $g$, and their indices are therefore upstairs: $g^{ij}$.
\end{itemize} 

\section{Strassen's Theorem over $\mathbb{R}^n$ revisited}

The aim of the present section is to introduce a version of Strassen's theorem and subsequently set the stage for a new fully constructive proof of it based on a notion of optimal kernels. Throughout, $\|\cdot\|$ will denote an arbitrary norm on $\mathbb{R}^n$. Moreover recall the following:

\textbf{Definition 1:} Given two Radon measures $\mu_1, \mu_2$ over $\mathbb{R}^n$ we shall say that $\mu_1 \preceq \mu_2$ if and only if
\begin{align}
\int_{\mathbb{R}^n}f(x) \, d\mu_1(x) \leq \int_{\mathbb{R}^n}f(x) \, d\mu_2(x)
\end{align}

for all convex functions $f : \mathbb{R}^n \rightarrow \mathbb{R}$. The partial order $\preceq$ is known as \textsl{convex order}.
\textbf{Definition 2:} Given an open set $\Omega \subset \mathbb{R}^n$, we define the space of continuous functions $\mathcal{E}(\Omega)$ with at most linear growth as
\begin{align}
\mathcal{E}(\Omega) := \K{\{}{\}}{f \in C(\Omega) \, | \, \sup_{x \in \Omega}\frac{|f(x)|}{1 + \|x\|} < \infty}.
\end{align}

\textbf{Theorem 1 (Strassen \cite{Strassen}):} Let $\mu_1, \mu_2 \in C_0(\mathbb{R}^n)^{*}$ Radon probability measures on $\mathbb{R}^n$, and let $\rho_1, \rho_2$ denote their generalized densities: $d\mu_1(x) := \rho_1(x)d^nx$, $d\mu_2(x) := \rho_2(x)d^nx$.
More strongly we require $\mu_1, \mu_2 \in \mathcal{E}(\mathbb{R}^n)^*$:
\begin{align}
\int_{\mathbb{R}^n}\|x\| d\mu_i(x) < \infty \s i = 1,2.
\end{align}

Then the following statements are equivalent:
\begin{itemize}
\item\textbf{Property 1:}
\begin{align}
\mu_1 \preceq \mu_2.
\end{align}

\item \textbf{Property 2:} There is a certain non-negative \textsl{kernel} $K(\, \cdot \, | \, \cdot \, ) : \mathbb{R}^n \times \mathbb{R}^n \rightarrow \mathbb{R}_{\geq 0}$ with the following properties:
\begin{align}
&\int_{\mathbb{R}^n} (1, y) \, K(y | x) \,  d^n y = (1, x) \s \textrm{a.s. w.r.t } \mu_1,\\
&\int_{\mathbb{R}^n} K(y | x)d\mu_1(x) = \rho_2(y).
\end{align}

More formally $K$ is a short-hand for measure-valued measures of the type
\begin{align}
\nu_{\bullet} : \mathcal{B}(\mathbb{R}^n) &\rightarrow \mathrm{RadonMeasures}(\mathbb{R}^n),\\
U &\mapsto \nu_U.
\end{align}

such that for all $U \in \mathcal{B}(\mathbb{R}^n)$
\begin{align}
\int_{\mathbb{R}^n} (1, y)d\nu_U(y) &= \int_U (1, x)d\mu_1(x),\\
\nu_{\mathbb{R}^n} &= \mu_2.
\end{align}
\end{itemize}

\medskip
\noindent
The proof clearly has two parts, $1 \Rightarrow 2$ and $2 \Rightarrow 1$. The second part is straightforward: 

\textbf{ Proof $2 \Rightarrow 1$  (Sketch) : } Let $f$ convex as in $1$, then by $2$ it follows
\begin{align}
&\int_{\mathbb{R}^n} f(y) d\mu_2(y)\\
&= \int_{\mathbb{R}^n \times \mathbb{R}^n} f(y) K(y | x) d\mu_1(x) d^n y\\
&\geq  \int_{\mathbb{R}^n} f\K{(}{)}{\int_{\mathbb{R}^n} y K(y | x) d^n y} d\mu_1(x)\\
&=  \int_{\mathbb{R}^n} f(x)d\mu_1(x).
\end{align}

More formally, choose a finite partition of $\mathbb{R}^n$ in $\mu_1$-measurable sets $U_{-1}, U_0, \dots, U_N$ such that $\mu_1(U_{-1}) = 0$, $\mu(U_k) > 0$ for $k = 0, \dots, N$, then
\begin{align}
\infty > &\int_{\mathbb{R}^n} f(y) d\mu_2(y)\\
&= \sum_{k = -1}^N \int_{\mathbb{R}^n} f(y) d\nu_{U_k}(y)\\
&\geq \sum_{k = 0}^N  f\K{(}{)}{\int_{\mathbb{R}^n} y \frac{d\nu_{U_k}(y)}{\mu_1(U_k)}} \mu_1(U_k)\\
&=\sum_{k = 0}^N  f\K{(}{)}{\frac{\int_{U_k} x \, d\mu_1(x)}{\mu_1(U_k)}} \mu_1(U_k)\\
&\rightarrow \int_{\mathbb{R}^n}  f(x) d\mu_1(x),
\end{align}

where the limit is taken over a a sequence of partitions. In particular one can choose convex sets with finite diameter and a set outside of a large compact set and use the continuity of convex functions.
\begin{align*}
\s\s\s\s\s\s\s\s\s\s\s\s\s\s\s\s\s\s\s\s\s\s\s\s\s\s\s\s\s\s\s\s\s\s\s\s\s\s\s\s\s\s\s\s\s\s\s\s\s\s\s\s\s\s\s\s\s\s\s\s\s\s\s\s\s\s\s \square
\end{align*}

For the harder $1 \Rightarrow 2$ implication, we will have to develop some machinery. The general idea is to construct a natural convex functional on the set of kernels satisfying \textsl{Property 2} and to show that such a set is not empty by establishing the existence of the minimum of such a functional. In the next section we construct the candidate kernel and in subsequent sections we show that such a kernel always exists if \textsl{Property 1} is fulfilled. The next section therefore can be viewed as a heuristic section and full rigour there is not attempted, rather it serves as an explanation of the form of the kernel we will choose for the subsequent proof. Throughout, \textsl{Property 1}, \textsl{Property 2} will refer to those of Theorem 1.

\section{Kernel from Mutual Information}
The functional we shall consider is the mutual information. In particular, given \textsl{Property 2}, $K$ is a conditional probability distribution, therefore if it exists, we can define the joint probability distribution
\begin{align}
d\mu_K(x,y) := K(y|x)d\mu_1(x)d^ny,
\end{align}

with marginals
\begin{align}
&d\mu_1(x) = \int_{y \in \mathbb{R}^n}K(y|x)d\mu_1(x)d^ny,\\
&d\mu_2(y) = \int_{x \in \mathbb{R}^n}K(y|x)d\mu_1(x)d^ny.\\
\end{align}

The mutual information is then given by
\begin{align}
S[K] &:= \int_{\mathbb{R}^n \times \mathbb{R}^n} d\mu_K(x,y) \log\K{(}{)}{\frac{d\mu_K(x,y)}{d\mu_1(x)d\mu_2(y)}}\\
&= \int_{\mathbb{R}^n \times \mathbb{R}^n}   \log\K{(}{)}{\frac{K(y|x)}{\rho_2(y)}} K(y|x)d\mu_1(x)d^ny,
\end{align}

which is non-negative and clearly convex in $K$. Recall that the mutual information expresses the maximum over all product distributions, of the likelyhood that a given joint distribution (when viewed as an empirical distribution) results (by sampling) from such product distributions. Such a maximum is attained at the product of the joint's marginals. Minimizing $S$ over $K$ means looking for the kernel that keeps the joint distribution as ``close as possible" to a product distribution. Hence if we do not impose any constraints on $K$, the minimum of $S$ is necessarily attained with $K(y|x) = \rho_2(y)$, where $S = 0$. We shall impose the linear constraints of \textsl{Property 2} by extending $S$ to include Lagrange multiplier functions $a \in L_{\mu_1}^1(\mathbb{R}^n), b_i \in L_{\mu^i_1}^1(\mathbb{R}^n),  c \in L_{\mu_2}^1(\mathbb{R}^n)$, where $d\mu_1^i(x) := x^i d\mu_1(x)$:
\begin{align}
\label{lossFirst}
S[K, a, b, c] &:= S[K] - \int_{\mathbb{R}^n} \K{\langle}{\rangle}{ (a(x), b(x)), \K{[}{]}{\int_{\mathbb{R}^n}(1, y)K(y|x)d^n y - (1, x)}}  d\mu_1(x)\\
&\s\s\s -\int_{\mathbb{R}^n} c(y) \K{[}{]}{\int_{\mathbb{R}^n} K(y|x)d\mu_1(x) - \rho_2(y)}d^n y.
\end{align}

The local minimum of $S$ is given by,
\begin{align}
\delta_KS[K^*, a, b, c] = 0
\end{align}

which is equivalent to:
\begin{align}
0 = \rho_1(x)\K{[}{]}{1 + \log\K{(}{)}{\frac{K^*(y|x)}{\rho_2(y)}} - a(x) - \langle b(x), y\rangle - c(y)}.
\end{align}

Therefore, if the local minimum exists, it is given by:
\begin{align}
K^*(y|x) = \rho_2(y) \exp\K{(}{)}{a(x) + \langle b(x), y\rangle + c(y) - 1}.
\end{align}

Moreover it is unique due to the strict convexity of the exponential as a function of $a,b,c$. We can now solve for the constraint that $K^*$ diffuses $\mu_1$ to $\mu_2$ to obtain an expression for $K^*$ that only depends on $a, b$:
\begin{align}
K^*(y|x) = \rho_2(y) \frac{\exp\K{(}{)}{a(x) + \langle b(x), y\rangle}}{\int_{\mathbb{R}^n} d\mu_1(u)\exp\K{(}{)}{a(u) + \langle b(u), y\rangle}}.
\end{align}

At this point Theorem 1 would follow from:

\textbf{Main Lemma: } Given \textsl{Property 1}, there is a sequence of kernels $K^*_m$ with corresponding $a^{(m)} \in L^{\infty}_{\mu_1}(\mathbb{R}^n, \mathbb{R}),  b^{(m)} \in  L^{\infty}_{\mu_1}(\mathbb{R}^n, \mathbb{R}^n)$ defining a sequence of measures $d\nu_{m,U}$, $U \in \mathcal{B}(U)$ that converges to a kernel $K$ as in \textsl{Property 2} in the following sense: There is a family of Radon measures $d\nu_U$ with $U \in \mathcal{B}(\mathbb{R}^n)$ such that
\begin{align}
\mathrm{a)} &\s \lim_{m \rightarrow \infty} \int_{\mathbb{R}^n} \, (1, y)d\nu_{m,U}(y)  =  \int_{U} (1, x) d\mu_1(x),\\
\mathrm{b)} &\s \lim_{m \rightarrow \infty} \int_{\mathbb{R}^n }  f(y)d\nu_{m,U}(y) =  \int_{\mathbb{R}^n } f(y) d\nu_U(y), \s \forall f \in \mathcal{E}(\mathbb{R}^n),\\
\mathrm{c)} &\s\s  \mu^f : \mathcal{B}(\mathbb{R}^n) \rightarrow [0, \infty)\\
&\s\s\s\s\s\s\s\s\s\s U \mapsto \int_{\mathbb{R}^n}f(x)d\nu_U(x)\\
&\s\s \mathrm{is} \; \mathrm{a} \; \mathrm{Radon} \; \mathrm{measure} \s \forall f \in \mathcal{E}(\mathbb{R}^n).
\end{align}

\section{Proof of Main Lemma}

In order to prove the main lemma we start by discretizing the problem. In particular we choose a finite partition of $\mathbb{R}^n$ into $\mu_1$-measurable sets $U_k$:
\begin{align}
\mathbb{R}^n = \bigsqcup_{k = -1}^N U_k
\end{align}

with $\mu_1(U_{-1}) = 0, \s \mu(U_k) > 0 \s \forall k = 0, \dots, N$. We define:
\begin{align}
p_k &:= \mu_1(U_k),\\
x_k &=  \frac{\int_{U_k} x \, d\mu_1(x)}{p_k} \s \forall k = 0, \dots, N.
\end{align}

We aim to show the finite version of the Main Lemma:

\textbf{Main Lemma (Finite Version) : } There is a sequence of simple functions $a^{(m)}, b^{(m)}$, as in the main lemma, subordinate to the above finite partition such that for all $k = -1, \dots, N$
\begin{align}
\label{mainlemmaid}
\lim_{m \rightarrow \infty}  \int_{\mathbb{R}^n} \, (1, y)d\nu_{m, U_k}^*(y)  = \int_{U_k} (1, x)d\mu_1(x).
\end{align}

\textbf{Proof:} W.l.o.g. we shall assume:
\begin{align}
\int_{\mathbb{R}^n} \, x \, d\mu_1(x) = 0 \in \mathbb{R}^n.
\end{align}

If (\ref{mainlemmaid}) held at finite $m$, then it would be equivalent to the existence of simple $a,b$ such that
\begin{align}
\int_{\mathbb{R}^n} (1, y) \, \frac{p_k \exp\K{(}{)}{a_k + \langle b_k, y \rangle}}{\sum_{l = 0}^N p_l \exp\K{(}{)}{a_l + \langle b_l, y\rangle}} \; d\mu_2(y) \, = (p_k, p_k \, x_k)\s \forall k = 0, \dots, N.
\label{maineq}
\end{align}

We notice that the kernel $K^{*}$ is invariant under a shift
\begin{align}
(a_k, b_k) \mapsto (a_k +  c, b_k + d),
\label{shiftsymmetry}
\end{align}

where $c \in \mathbb{R}, d \in \mathbb{R}^n$. Therefore, without loss of generality we can choose $(a_k, b_k) = 0$ for some arbitrary $k$. For convenience we will choose:
\begin{align}
(a_0, b_0) = (0, 0).
\end{align}

Moreover, without loss of generality, since we can enumerate $U_k$ as we like, we can assume that
\begin{align}
0 = a_0 \geq a_k \s \forall k = 1, \dots, N.
\end{align}

We also notice that it suffices to show (\ref{maineq}) for $k = 1, \dots, N$, as the identity (at finite $m$) for $k = 0$ then follows from
\begin{align}
\sum_{k = 0}^N p_k (1, x_k) &= \int_{\mathbb{R}^n} (1, x)d\mu_1(x)\\
&= \int_{\mathbb{R}^n} (1, y) d\mu_2(y)\\
&= \sum_{k = 0}^N \int_{\mathbb{R}^n} (1, y) \, \frac{p_k \, \exp\K{(}{)}{a_k + \langle b_k, y \rangle}}{\sum_{l = 0}^N p_l \exp\K{(}{)}{a_l + \langle b_l, y\rangle}}\;d\mu_2(y),
\end{align}

where in the second identity we have used \textsl{Property 1} applied to the functions $f = 1, x_1, ..., x_n$ that are all both convex and concave, therefore the inequality for arbitrary convex functions becomes an equality.
Hence we are left to show (in the limit of $m\rightarrow \infty$):
\begin{align}
\int_{\mathbb{R}^n} (1, y) \, \frac{p_k \exp\K{(}{)}{a_k + \langle b_k, y \rangle}}{\sum_{l = 0}^N p_l \exp\K{(}{)}{a_l + \langle b_l, y\rangle}}\;d\mu_2(y)  = (p_k, p_k x_k) \s \forall k = 1, \dots, N.
\label{maineq_2}
\end{align}

To this aim we define the function:
\begin{align}
F : \mathbb{R}^{N(n + 1)} &\rightarrow \mathbb{R}^{N(n+1)}\\
F((a_1, b_1), \dots, (a_N, b_N))_k &:=  \int_{\mathbb{R}^n} (1, y) \, \frac{p_k \exp\K{(}{)}{a_k + \langle b_k, y \rangle}}{\sum_{l = 0}^N p_l \exp\K{(}{)}{a_l + \langle b_l, y\rangle}}\;d\mu_2(y)\\
&= \nabla \Phi((a_1, b_1), \dots, (a_N, b_N)),
\end{align}

where:
\begin{align}
\Phi((a_1, b_1), \dots, (a_N, b_N)) := \int_{\mathbb{R}^n}  \log\K{(}{)}{ p_0 + \sum_{l = 1}^N p_l \exp\K{(}{)}{a_l + \langle b_l, y\rangle}} \; d\mu_2(y).
\end{align}

Now, without loss of generality, we can assume that $\rho_2$ is not concentrated on any codimension $1$ hyperplane $H \subset \mathbb{R}^n$. Meaning for all such $H$ $\mu_2(H) < 1$. If this were not the case then so too would $\mu_1(H) = 1$, otherwise we could construct a convex function that violates \textsl{Property 1}. Namely, we could choose the function:
\begin{align}
f(x) := \max\K{(}{)}{\langle x - x_H, n_H\rangle, 0},
\end{align}

where $x_H$ is an arbitrary point on $H$ and $n_H$ is the normal versor to $H$ in either orientation. Then the main lemma would reduce altogether to the lower dimensional setting $n \mapsto n - 1$. Hence, in conclusion we shall assume in what follows that:
\begin{align}
\mu_2(H) < 1 \s \forall \textrm{ hyperplanes } H.
\end{align}

Then it follows that $\Phi$ is strictly convex. To prove this, we shall choose an arbitrary vector $(v_1, \dots, v_N) \in \mathbb{R}^{N(n+1)}$ with $v_k \in \mathbb{R}^{n+1}$ and evaluate the Hessian of $\Phi$ as a quadratic form on it:
\begin{align}
\langle v, H(\Phi)((a,b))v\rangle = \int_{\mathbb{R}^n}\K{[}{]}{\mathbb{E}_y \K{[}{]}{\langle v_{\bullet}, \hat{y}\rangle^2} - \mathbb{E}_y \K{[}{]}{\langle v_{\bullet}, \hat{y}\rangle}\mathbb{E}_y \K{[}{]}{\langle v_{\bullet}, \hat{y}\rangle}} \; d\mu_2(y) ,
\end{align}

where $\hat{y} := (1,y)$, $v_0 = 0$ and
\begin{align}
\mathbb{E}_y \K{[}{]}{\xi_{\bullet}} := \sum_{l = 0}^N \xi_l p_l(y) :=  \frac{\sum_{l = 0}^N  \xi_l  p_l\exp\K{(}{)}{a_l + \langle b_l, y\rangle}}{\sum_{l = 0}^N p_l \exp\K{(}{)}{a_l + \langle b_l, y\rangle}}.
\end{align}

Without loss of generality $N > 0$ (otherwise $\mu_2$ would be a Dirac measure, and hence so would $\mu_1$) and since $p_k > 0$ for all $k = 0, \dots, N$, and $v_0 = 0$, it follows
\begin{align}
&\langle v, H(\Phi)((a,b))v\rangle = 0 \Leftrightarrow \\
&\exists v \in \mathbb{R}^{N(n+1)}, y_0 \in \mathbb{R}^n \s s.t. \s  \langle v_k, \hat{y} - \hat{y}_0\rangle = 0 \s \forall k = 1, \dots, N, y \in \textrm{supp}(\mu_2).
\end{align}

However the r.h.s is true iff there is a hyperplane $H \subset \mathbb{R}^n$ such that $\mu_2(H) = 1$, which we have ruled out. Hence, finally, w.l.o.g.\ we can assume that $\Phi$ is strictly convex. Now we shall use the fact that for a strictly convex function $\Phi : \mathbb{R}^d \rightarrow \mathbb{R}$ (in this case $d = N(n+1)$), the function
\begin{align}
F = \nabla \Phi
\end{align}

is a diffeomorphism onto its image. From the inverse function theorem, it suffices to prove that $F$ is injective and that its Jacobian is invertible everywhere. The last part we just showed, so it remains to show that $F$ is injective: Since $\Phi$ is strictly convex, $\forall x \neq y \in \mathbb{R}^n$:
\begin{align}
\Phi(y) - \Phi(x) &> \lim_{\lambda \rightarrow 0} \frac{\Phi((1 - \lambda)x + \lambda y) - \Phi(x)}{\lambda}\\
&= \langle y - x, \nabla \Phi(x) \rangle.
\end{align}

Hence also interchanged
\begin{align}
\Phi(x) - \Phi(y) > \langle x - y, \nabla \Phi(y) \rangle.
\end{align}

Finally
\begin{align}
0 < \langle y - x, F(y) - F(x)\rangle.
\end{align}

Hence $x \neq y$ implies $F(x) \neq F(y)$.

Now that we know that $F$ is a diffeomorphism onto its image and we also know that in particular $F$ is an open map in the neighborhood of $0 \in \mathbb{R}^{N(n+1)}$ where
\begin{align}
F(0) = ((p_1, 0), (p_2, 0), \dots, (p_N, 0)),
\end{align}

which means that if we consider the open half-line
\begin{align}
\mathcal{L} := \{X_{\mu} := ((p_1, \mu p_1x_1), (p_2, \mu p_2x_2), \dots, (p_N, \mu p_N x_N)) \, | \, \mu > 0\},
\end{align}

then
\begin{align}
\mathcal{L} \cap F(\mathbb{R}^{N(n + 1) }) \neq \emptyset.
\end{align}

Since $\overline{F(\mathbb{R}^{N(n + 1) - 1})}$ is compact ( follows easily straight from the definition and the fact that $\mu_2$ has finite mean ), it follows that:
\begin{align}
\mathcal{L} \cap \partial F(\mathbb{R}^{N(n + 1)} ) \neq \emptyset.
\end{align}

Moreover, since $F$ is a diffeomorphism:
\begin{align}
\partial F(\mathbb{R}^{N(n + 1) }) = F(S_{\infty}^{N(n + 1) -1}),
\end{align}

where $S_{\infty}^{N(n + 1) -1}$ is the sphere at infinity centered at $0$. At this point, the claim of the proposition follows if we can show:
\begin{align}
X_{\mu} \in \mathcal{L} \cap F(S_{\infty}^{N(n + 1) -1}) \Leftrightarrow \mu \geq 1,
\end{align}

because this would then imply:
\begin{align}
X_1 = ((p_1, p_1x_1), \dots, (p_N, p_Nx_N)) \in \overline{F(\mathbb{R}^{N(n + 1) })},
\end{align}

which is our claim. A point in $F(S_{\infty}^{N(n + 1)-1})$ is of the form
\begin{align}
\lim_{r\rightarrow \infty} F((ra_1, rb_1), \dots, (r a_N, r b_N))_k = \int_{\mathbb{R}^n}  (1, y) \prod_{l \neq k} 1_{a_k - a_l + \langle b_k - b_l, y \rangle \geq 0} \;d\mu_2(y),
\end{align}

where $((a_1, b_1), \dots, (a_N, b_N)) \neq 0$. Therefore, a point $X_{\mu} \in \mathcal{L} \cap F(S_{\infty}^{N(n + 1)-1})$ must satisfy
\begin{align}
p_k(1, \mu x_k) = \int_{\mathbb{R}^n} (1, y) \prod_{l \neq k} 1_{a_k - a_l + \langle b_k - b_l, y \rangle \geq 0}\; d\mu_2(y).
\end{align}

Therefore, summing up and recalling that $(a_0, b_0) = 0$,
\begin{align}
&\sum_{k = 0}^N p_k \K{(}{)}{a_k + \mu \langle b_k, x_k \rangle} \\
&= \int_{\mathbb{R}^n}  \sum_{k = 0}^N(a_k + \langle b_k, y \rangle) \prod_{l \neq k} 1_{a_k - a_l + \langle b_k - b_l, y \rangle \geq 0}\;d\mu_2(y)\\
&= \int_{\mathbb{R}^n}  \max\K{(}{)}{0, a_1 + \langle b_1, y\rangle, \dots, a_N + \langle b_N, y\rangle}d\mu_2(y)\\
&\geq  \int_{\mathbb{R}^n}\max\K{(}{)}{0, a_1 + \langle b_1, y\rangle, \dots, a_N + \langle b_N, y\rangle} d\mu_1(y) \\
&\geq \sum_k p_k \max\K{(}{)}{0, a_1 + \langle b_1, x_k\rangle, \dots, a_N + \langle b_N, x_k\rangle},
\end{align}

where we have used \textsl{Property 1} and the convexity of the max function twice. It follows that:
\begin{align}
\mu \sum_{k = 0}^N p_k\langle b_k, x_k\rangle &\geq \sum_k p_k \K{(}{)}{\max\K{(}{)}{0, a_1 + \langle b_1, x_k\rangle, \dots, a_N + \langle b_N, x_k\rangle} - a_k}.
\label{final_inequality}
\end{align}

Since we could choose $0 = a_0 \geq a_k$ for all $k$:
\begin{align}
\mu \sum_{k = 0}^N p_k\langle b_k, x_k\rangle &\geq \sum_{k = 0}^N p_k\langle b_k, x_k\rangle,
\end{align}

and
\begin{align}
\mu \sum_{k = 0}^N p_k\langle b_k, x_k\rangle &\geq 0.
\end{align}

It follows that either $\sum_{k = 0}^N p_k\langle b_k, x_k\rangle > 0$ in which case $\mu \geq 1$ and the proposition is proved, or $\sum_{k = 0}^N p_k\langle b_k, x_k\rangle = 0$. However in that latter case, it would result from (\ref{final_inequality}) that $a_k = 0$ for all $k = 0, \dots, N$ hence it would follow from (\ref{final_inequality}) that
\begin{align}
0 = \int_{\mathbb{R}^n} \max\K{(}{)}{0, \langle b_1, y\rangle, \dots, \langle b_N, y\rangle}d\mu_2(y).
\end{align}

Given the hyperplane assumption on $\mu_2$, such can only happen if $b_1 = \dots = b_N = 0$ but in that case $((a_1, b_1), \dots, (a_N, b_N)) = 0$ which contradicts the assumption $X_{\mu} \in F(S_{\infty}^{N(n+1)-1})$.
Moreover if \textsl{Property 1} holds strictly on convex functions that are non-affine then $\mu > 1$, because one such function is the max function employed above which would give a strict inequality. In that case $a,b$ are bounded.
\begin{align*}
\s\s\s\s\s\s\s\s\s\s\s\s\s\s\s\s\s\s\s\s\s\s\s\s\s\s\s\s\s\s\s\s\s\s\s\s\s\s\s\s\s\s\s\s\s\s\s\s\s\s\s\s\s\s\s\s\s\s\s\s\s\s\s\s\s\s\s \square
\end{align*}

\textbf{End of Proof of Main Lemma, hence of Theorem 1.}

\textbf{Proof of a):}
We choose a sequence of increasing concentric cubes centered around $0$, $C_N := [-2^{N-1}, 2^{N-1}]^n$. We then partition $C_N$ into $P_N := 4^{Nn}$ cubes $V_{N,k}$ , $k = 1, \dots, P_N$ each of volume $2^{-Nn}$. By Main Lemma (Finite Version), for all $N \geq 0, \epsilon(N) > 0$, subordinate to the partition of $\mathbb{R}^n$ defined by  $\mathbb{R}^n \backslash C_N$ and $V_{N, k}$, we have a sequence $a_{N, m}, b_{N, m}$ of simple functions such that for all $V_{N, k}$, $k = 1, \dots, P_N$, there is $M > 0$, such that for all $m > M$
\begin{align}
\K{\|}{\|}{\int_{\mathbb{R}^n}(1, y)d\nu^*_{m, V_{N, k}}(y) - \int_{V_{N,k}} (1, x)d\mu_1(x)} < \epsilon(N).
\end{align}

For $U \in \mathcal{B}(\mathbb{R}^n)$, we consider all cubes $V_{k,N}$ that are contained in $U$ and define $K_N(U)$ to be the compact set given by their union. Then:
\begin{align}
&\K{\|}{\|}{\int_{\mathbb{R}^n}(1, y)d\nu^*_{m, K_{N}(U)}(y) - \int_{U} (1, x)d\mu_1(x)} \\
&\leq \K{\|}{\|}{\int_{\mathbb{R}^n}(1, y)d\nu^*_{m, K_{N}(U)}(y) - \int_{K_N(U)} (1, x)d\mu_1(x)} + \K{\|}{\|}{\int_{K_N(U)} (1, x)d\mu_1(x) - \int_{U} (1, x)d\mu_1(x)}\\
&\leq Z_N(U) \epsilon(N) + \K{\|}{\|}{\int_{U \backslash K_N(U)} (1, x)d\mu_1(x) },
\end{align}

where $Z_N(U)$ is the number of $V_{N, k}$ cubes contained in $U$. Since the measure corresponding to $(1 + \|x\|) d\mu_1(x)$ is a finite Radon measure, it is inner regular on Borel sets, hence the second term can be chosen arbitrarily small with increasing $N$, and so can the first. In particular we can extract a diagonal subsequence
\begin{align}
\nu_{N, U} := \nu^*_{m_N, K_N(U)}
\end{align}

such that
\begin{align}
\lim_{N \rightarrow \infty} \int_{\mathbb{R}^n}(1, y)d\nu_{N,U} = \int_{U} (1, x)d\mu_1(x).
\end{align}

\textbf{Proof of b):} We observe that $d\nu_{N,U}$  is of the form
\begin{align}
d\nu_{N,U}(y) = d\mu_2(y)\xi_{N,U}(y),
\end{align}

where $\xi_{N,U}(y) \in [0,1]$ for all $y \in \mathbb{R}^n$. Therefore the sequence of measures $\nu_{N,U}(y)$ is uniformly bounded on $\mathcal{E}(\mathbb{R}^n)$ and tight. Hence by Prokhorov's theorem, it has a convergent subsequence in the weak-$^{\star}$ topology. Without loss of generality, we can replace the sequence with such a convergent subsequence. In other words, there is a Radon measure $\nu_U \in \mathcal{E}(\mathbb{R}^n)^*$ with 
\begin{align}
\lim_{N \rightarrow \infty} \int_{\mathbb{R}^n} f(y) d\nu_{N, U}(y) = \int_{\mathbb{R}^n }  f(y) d\nu_U(y), \s \forall f \in \mathcal{E}(\mathbb{R}^n).
\end{align}

\textbf{Proof of c):} We denote by $\mathcal{A}_N$ the algebra generated by $V_{N,k}$, $k = 1, \dots, P_N$. It holds $\mathcal{A}_N \subset \mathcal{A}_{N+1}$ and the union $\mathcal{A} := \bigcup_{N \geq 0}\mathcal{A}_N$ generates the Borel algebra $\mathcal{B}(\mathbb{R}^n)$. Moreover $\mathcal{A}$ is countable, and we shall enumerate the elements of $\mathcal{A}$ as $U_i$, $i \in \mathbb{N}$ such that if $U_i \in \mathcal{A}_N$, $\emptyset \neq U_j \in \mathcal{A}_{N+1}\backslash \mathcal{A}_N$, then $i < j$. From part b) we know that for each $U_i$ there is a subsequence $\nu^*_{m_{N_i}, K_{N_i}(U_i)}$ converging weak-$^{\star}$ to a $\nu_{U_i}$. We now construct an increasing sequence $n_k \in \mathbb{N}$ such that $\nu^*_{m_{n_k}, K_{n_k}(U_i)}$ converges for all $U_i$. The latter is achieved by defining iteratively $n_k^{1}$ such that the subsequence $\nu^*_{m_{n_k^1}, K_{n_k^1}(U_1)}$ converges, $n_k^{2}$ a subsequence of $n_k^{1}$ such that $\nu^*_{m_{n_k^2}, K_{n_k^2}(U_2)}$ converges and so on, finally $n_k := n_k^k$. Then for all $N \geq 0$, there is $k_N > 0$ monotonic in $N$ such that for all $k \geq k_N$,
\begin{align}
K_{n_k}(U) = U, \s \forall U \in \mathcal{A}_N
\end{align}

and for all $f \in \mathcal{E}(\mathbb{R}^n)$, the map:
\begin{align}
\widetilde{\mu}_{k}^f : \mathcal{A}_N &\rightarrow [0, \infty)\\
U &\mapsto \int_{\mathbb{R}^n} f(x) d\nu_{m_{n_k}, U}(x)
\end{align}

is uniformly bounded, as a result of
\begin{align}
\widetilde{\mu}_k^f(\mathbb{R}^n) = \int_{\mathbb{R}^n}f(x)d\mu_2(x),
\end{align}

defines a pre-measure on $\mathcal{A}_N$ and for all $l_N \geq k_N$,
\begin{align}
\label{limittight}
\lim_{N \rightarrow \infty} \widetilde{\mu}^f_{l_N}(\mathbb{R}^n \backslash C_N) \rightarrow 0.
\end{align}

Indeed, let $f$ such that $f/(1 +\|\cdot\|)$ has unit sup norm. Then for all $R>0, \epsilon >0$, there is  $M \geq 0$ such that for all $N \geq M$
\begin{align}
\widetilde{\mu}^f_{l_N}(\mathbb{R}^n \backslash C_N) &\leq \epsilon + \int_{\mathbb{R}^n \backslash B_R(0)} (1 + \|y\|)d\mu_2(y)+ (1 + R)\int_{\mathbb{R}^n \backslash C_N} d\mu_1(x).
\end{align}

Hence the limit (\ref{limittight}) follows. As a result, the limit of $\nu^*_{m_{n_k}, K_{n_k}(\cdot)}$ defines a pre-measure $\widetilde{\mu}^f$ on $\mathcal{A}$ for all $f \in \mathcal{E}(\mathbb{R}^n)$. Since $\widetilde{\mu}^f$ is finite, from the Carath\'eodory extension theorem, it extends to a unique measure $\mu^f$ on $\mathcal{B}(\mathbb{R}^n)$. Moreover by construction $\mu^f$ is Radon.
\begin{align*}
\s\s\s\s\s\s\s\s\s\s\s\s\s\s\s\s\s\s\s\s\s\s\s\s\s\s\s\s\s\s\s\s\s\s\s\s\s\s\s\s\s\s\s\s\s\s\s\s\s\s\s\s\s\s\s\s\s\s\s\s\s\s\s\s\s\s\s \square
\end{align*}



\section{Hessian Martingales}\label{hessianmartingales}
In this section we will specialize to kernels obtained from the mutual information that are smooth, and completely characterize them. We first start with the archetypal example of such a kernel, which is that defined by Brownian Motion. Indeed, consider Brownian Motion in dimension $1$. In that case the kernel between times $s < t$ is given by:
\begin{align}
K(y|x) &= \frac{1}{\sqrt{2\pi (t - s)}}\exp\K{(}{)}{-\frac{(y - x)^2}{2(t - s)}},
\end{align}

which is of the form:
\begin{align}
\exp\K{(}{)}{a(x) + b(x)y + c(y)}.
\end{align}

Therefore we can reinterpret the Brownian kernel as the unique martingale kernel diffusing the marginal $\rho_s$ to $\rho_t$, which minimizes the mutual information of the joint distribution, namely it keeps it ``as close as possible" to a product distribution. Similarly Brownian motion in arbitrary dimension with a general constant covariance matrix is also optimal w.r.t. mutual information. More generally, we can look at optimal kernels on (an open simply-connected set with sufficiently regular boundary of) $\mathbb{R}^n$:
\begin{align}
K(y|x) = \exp\K{(}{)}{a(x) + \langle b(x), y\rangle + c(y)},
\end{align}

where $a,b$ are smooth over all of $\mathbb{R}^n$. Then, it follows from the constraints:
\begin{align}
0 &= \partial_{k} \int_{\mathbb{R}^n} \exp\K{(}{)}{a(x) + \langle b(x), y\rangle + c(y)}d^n y\\
&= \partial_{k} a(x) +  \langle \partial_k b(x), x\rangle,
\end{align}

and
\begin{align}
\delta_k^l &= \partial_{k} \int_{\mathbb{R}^n} y^l \exp\K{(}{)}{a(x) + \langle b(x), y\rangle + c(y)}d^n y\\
&= x^l \partial_{k} a(x) +  \int_{\mathbb{R}^n}\langle \partial_k b(x), y\rangle y^l \exp\K{(}{)}{a(x) + \langle b(x), y\rangle + c(y)}d^n y.
\end{align}

From the first constraint it follows
\begin{align}
b_k(x) = \partial_k \K{(}{)}{ a(x) + \langle b(x), x\rangle }.
\end{align}

Therefore, there is $\phi$ such that
\begin{align}
b = \nabla \phi,
\end{align}

then, w.l.o.g. (absorbing the constant term in $a$)
\begin{align}
a(x) = \phi(x) - \langle \nabla \phi(x), x\rangle.
\end{align}

The second, martingality constraint then becomes
\begin{align}
\delta_k^l &= - x^r x^l \partial_k \partial_r \phi(x) + \partial_k\partial_r \phi(x) \int_{\mathbb{R}^n} y^r y^l K(y|x)d^n y\\
&= \partial_k\partial_r \phi(x) \mathrm{Cov}_{K}(y^r, y^l \, | \, x).
\end{align}

Hence the Hessian $g(x) := H(\phi)(x)$ is positive definite for all $x \in \mathbb{R}^n$ and
\begin{align}
g^{-1}(x) = \mathrm{Cov}_{K}(y, y^T \, | \, x).
\label{conditionalcovariance}
\end{align}

Hence finally there is a smooth strictly convex function $\phi$ such that
\begin{align}
K(y|x) = \exp\K{(}{)}{\phi(x) + \langle \nabla \phi(x), y - x\rangle + c(y)}.
\end{align}

We will now show how $\phi(x)$ is uniquely determined by $c(y)$ and as a byproduct we determine the precise class of allowed strictly convex functions $\phi$. We start by recalling that, by the shift symmetry of the mutual information with constraints ($(a(x), b(x)) \mapsto (a(x) + a_c, b(x) + b_c) \s a_c \in \mathbb{R}, b_c \in \mathbb{R}^n$) we can choose $x_*$ such that
\begin{align}
&\phi(x_*) = 0,\\
&\nabla \phi(x_*) = 0.
\end{align}

In practice $\phi$ will be bounded within an open convex subset $\Omega \subset \mathbb{R}^n$ and will diverge to $+\infty$ on $\partial \Omega$. Therefore we can choose $x_* \in \Omega$. A natural choice is therefore
\begin{align}
x_* = \int_{\mathbb{R}^n} \, x \, d\mu_1(x).
\end{align}

It follows
\begin{align}
d\mu(y|x_*) := K(y | x_*)d^n y = \exp\K{(}{)}{c(y)}d^n y .
\end{align}

Hence
\begin{align}
\int_{\mathbb{R}^n} \exp\K{(}{)}{\langle \nabla \phi(x), y \rangle} d\mu(y \, | \, x_*) = \exp\K{(}{)}{\langle x, \nabla \phi(x)\rangle - \phi(x)}.
\end{align}

Since $\phi$ is strictly convex in $\Omega$, the map
\begin{align}
F := \nabla \phi : \Omega \rightarrow \mathbb{R}^n
\end{align}

is a diffeomorphism onto its image\footnote{See e.g. the proof of the Main Lemma.}. Let $\Omega^* := F(\Omega)$. Then it follows
\begin{align}
\int_{\mathbb{R}^n}  \exp\K{(}{)}{\langle k, y \rangle}d\mu(y \, | \, x_*)  &= \exp\K{(}{)}{\langle k, F^{-1}(k)\rangle - \phi(F^{-1}(k))}\\
&= \exp\K{(}{)}{\psi(k)} \s \forall k \in \Omega^*,
\end{align}

where $\psi : \Omega^* \rightarrow \mathbb{R}$ is the Legendre transform of $\phi$
\begin{align}
\psi(k) := \sup_{x \in \Omega} \K{(}{)}{\langle k, x\rangle - \phi(x)}.
\end{align}

Hence in particular $\psi$ is convex. On the other hand for any (Radon) probability measure $d\nu(y)$ over $\mathbb{R}^n$, the function:
\begin{align}
\psi(k) := \log\K{(}{)}{\int_{\mathbb{R}^n}\exp\K{(}{)}{\langle k, y\rangle}d\nu(y)}
\label{psi}
\end{align} 

is convex on its domain of definition $\Omega^*$ defined by $k \in \mathbb{R}^n$ where $|\psi(k)| < \infty$. Therefore, rather than starting from $\phi$ strictly convex, we can turn the problem upside down and define the kernel completely from a choice of $d\nu(y)$. Then $\phi$ is obtained as the Legendre transform of $\psi$:
\begin{align}
\phi(x) := \sup_{k \in \Omega^*}\K{(}{)}{\langle x, k\rangle - \psi(k)}.
\end{align}

\textbf{Remark: }Note that not all (strictly) convex functions are of type (\ref{psi}). The latter are such that $\exp\K{(}{)}{\psi(k)}$ admits an analytic continuation $\exp\K{(}{)}{\psi(ik)}$ for all $k \in \mathbb{R}^n$ as that is indeed the characteristic function of $\nu$. In particular $\mathrm{Re}(\psi(ik))$ is well defined (independent of the branch cut of the log) and $\mathrm{Re}(\psi(ik)) \leq 0$ for all $k \in \mathbb{R}^n$. Instead, e.g. $\psi(k) := k^2 + k^4$ on $\mathbb{R}$ is strictly convex, but $\mathrm{Re}(\psi(ik)) = -k^2 + k^4 > 0$ for $k > 1$.

Finally we will refer to the optimal kernels as \textsl{discrete time Hessian (local) martingales}.

\section{Hessian Martingales in Continuous Time}

The present section does not attempt to be rigorous, but rather to briefly sketch a conjectural structure of Hessian martingales in continuous time. We shall consider a one-parameter family of marginals $\{\mu_t\}_{t \in [0, T)}$ with $T > 0$ such that $\mu_s \preceq \mu_t$ for $s \leq t$\footnote{Such a family is also known as a ``peacock" (processus croissant pour l' \`ordre convexe).}. We shall call such a family differentiable if for any pair of times $s \leq  t \in [0, T)$ the corresponding $\phi \in C^2(\Omega)$ on its domain $\Omega$. For each $t, t + \epsilon \in [0, T)$ with $\epsilon > 0$, we then have functions $\phi_{t, \epsilon}$, $\psi_{t, \epsilon}$ and measure $d\nu_{t, \epsilon}$ relative to the diffusion from $t$ to $t + \epsilon$ with:
\begin{align}
\psi_{t, \epsilon}(k) = \epsilon \log \K{(}{)}{\int_{\mathbb{R}^n}\exp\K{(}{)}{\frac{1}{\epsilon}\langle k, y\rangle}d\nu_{t, \epsilon}(y)}.
\end{align}

We shall define:
\begin{align}
&\phi_t(x) := \lim_{\epsilon \rightarrow 0} \sup_{k \in \Omega_{t, \epsilon}^*}\K{(}{)}{\langle k, x\rangle - \psi_{t, \epsilon}(k)}.
\end{align}

\textbf{Definition Attempt:} we shall say that $\{d\mu_t\}_{t \in [0, T)}$ admits a \textsl{continuous Hessian (local) martingale} process if $\phi_t \in C^2(\Omega)$, everywhere strongly convex, where $\Omega$ is an open convex set of $\mathbb{R}^n$ and $\K{.}{|}{\phi_t}_{\partial \Omega} = \infty$ for all $t \in [0,T)$. We have the following conjecture:

\textbf{Conjecture:} $\{d\mu_t\}_{t \in [0, T)}$  admits a continuous Hessian martingale process if and only if the densities $\{\rho_t\}_{t \in [0, T)}$  satisfy the following forward Kolmogorov equation:

\begin{align}
\label{forwardKolmogorovHessian}
\partial_t \rho_t(x) = \frac{1}{2}\partial_{i}\partial_j \K{(}{)}{g_t^{ij}(x) \rho_t(x)},
\end{align}

where
\begin{align}
\label{HessianMetricContinuous}
(g_t)_{ij}(x) = \partial_i\partial_j \phi_t(x),
\end{align}

with $\phi_t$ as specified in the ``definition attempt". We leave a proper formulation and proof to future work, however we would like to mention the intuition behind it. Namely that such a statement should be a consequence of (\ref{conditionalcovariance}) coupled with the Central Limit Theorem.

We shall formally express:
\begin{align}
K_{t, \epsilon}(y|x_*) =: \exp\K{(}{)}{-\frac{1}{\epsilon}S_{t, \epsilon}(y)}.
\end{align}

Then it follows from The Kolmogorov equation that
\begin{align}
\partial_{\epsilon}K_{t, \epsilon}(y|x_*) = \frac{1}{2}\partial_{i}\partial_j \K{(}{)}{g_{t + \epsilon}^{ij}(y) K_{t, \epsilon}(y|x_*)},
\end{align}

which we can interpret as a ``renormalization group" equation for $S_{t, \epsilon}$. Finally we shall remark that formally, equations (\ref{forwardKolmogorovHessian}), (\ref{HessianMetricContinuous}) result from the minimization of the continuum version of (\ref{lossFirst}):
\begin{align}
S_t[g_t^{-1}, \phi_t] :=&-\frac{1}{2}\int_{\mathbb{R}^n} \log \det \K{(}{)}{g_t^{-1}(x)}\rho_t(x)d^n x\\
& - \int_{\mathbb{R}^n} \phi_t(x) \K{(}{)}{\partial_t \rho_t(x) - \frac{1}{2}\partial_i \partial_j \K{(}{)}{g_t^{ij}(x) \rho_t(x)}}d^n x,
\end{align}

where here $\phi_t$ is viewed as a Lagrange multiplier. In words, the optimal (most unbiased) covariance matrix is the one that tends to be ``as non-degenerate as possible" relative to $\rho_t$ thus tending to smoothen the resulting process.
\section{Hessian Martingales from Incomplete Marginals}\label{IncompleteMarginals}

In this section we revisit the construction of (discrete time) Hessian martingales for the case when the marginals are not fully known. Instead, one has knowledge of the expectation value of a finite number of continuous functions in $\mathcal{E}(\mathbb{R}^n)$ for each time point. From the mathematical finance perspective, such functions can be regarded as (discounted) payoffs of European contingent claims, and their expectations correspond to risk neutral prices. The main results of this section will be \textsl{Theorem 2} and \textsl{Theorem 3}. The former will establish sufficient conditions for the existence of a Hessian martingale that can be constructed sequentially from early to later times. In particular Theorem 2 will show that there is an obstruction to such sequential construction, thus invalidating a sequential approach in practice\footnote{In fact, the aforementioned obstruction is generic for $n \geq 2$, but it is absent for standard practical cases when $n=1$.}. Theorem 3, instead will establish the existence and uniqueness of a Hessian martingale and detail its general global (in particular non-sequential) structure. The latter forms the basis of the practical approach to calibrating Hessian martingales in general. Before formulating Theorems 2,3, we recall the following instrumental concepts and results concerning convex functions.

\textbf{Definition 3:} Given a function $f : \Omega \rightarrow \mathbb{R}$, with $\Omega \subset \mathbb{R}^n$ open,  its \textsl{lower convex envelope} is given by:
\begin{align}
\mathrm{conv}_{\Omega}(f) : \Omega &\rightarrow \mathbb{R}\\
x &\mapsto \sup \{a + \langle b,x \rangle \, | \, a + \langle b,y \rangle \leq f(y) \, \forall y \in \Omega, a \in \mathbb{R}, b \in \mathbb{R}^n \}.
\end{align}

That is $\mathrm{conv}_{\Omega}(f)$ is the supremum over all convex functions majorized by $f$ in the set $\Omega$.

\textbf{Proposition 1:} $\mathrm{conv}_{\Omega}(f)$ is the Legendre transform of
\begin{align}
\psi(k) := \sup_{y \in \Omega}\K{(}{)}{\langle k, y\rangle - f(y)}.
\end{align}

\textbf{Proof:} see appendix \ref{convexenvelopes}. 

\textbf{Definition 4:} We say that a vector space $\mathcal{H}$ of functions $f : \mathbb{R}^n \rightarrow \mathbb{R}$ is \textsl{closed under lower convex envelopes relative to $\Omega$} if for all $f \in \mathcal{H}$ then either
\begin{itemize}
\item $\mathrm{conv}_{\Omega}(f) \in \mathcal{H}$ \s or
\item $\mathrm{conv}_{\Omega}(f) = -\infty$.
\end{itemize}

We will also need the following:

\textbf{Proposition 2:} Let $\Omega \subset \mathbb{R}^n$ open, $h\in C(\Omega)$, $\rho \in L^1(\mathbb{R}^n)$ a probability density with $\mathrm{supp}(\rho) = \Omega$ and $\phi_h$ the Legendre transform of
\begin{align}
\psi(k) := \log\K{(}{)}{\int_{\mathbb{R}^n}\rho(y)\exp\K{(}{)}{\langle k, y\rangle - h(y)}d^n y}.
\end{align}

Then for all $x \in \langle \Omega\rangle$ (the convex hull of $\Omega$)
\begin{align}
\lim_{r \rightarrow \infty} \frac{1}{r}\phi_{rh}(x) = \mathrm{conv}_{\Omega}(h)(x),
\end{align}
and $\phi_{h}, \mathrm{conv}_{\Omega}(h) \in C(\langle \Omega\rangle) \cup \{-\infty\}$. Moreover, if $h \in \mathcal{E}(\Omega)$, then $\mathrm{conv}_{\Omega}(h) \in \mathcal{E}(\Omega) \cup \{-\infty\}$.

\textbf{Proof:} see appendix \ref{convexenvelopes}.

\subsection{Sequential Calibration with Incomplete Marginals}

We start by considering two marginals $\mu_1 \preceq \mu_2$ that are not fully known. Instead what is known are $\mu_1$, while concerning $\mu_2$ one has the expectation value, w.r.t. it, of a class of continuous functions $h \in \mathcal{H} \subset \mathcal{E}(\mathbb{R}^n)$ including the affine functions. That is the following constraints are satisfied:
\begin{align}
\int_{\mathbb{R}^n}h(y)d\mu_2(y) = F[h] \s \forall h \in \mathcal{H},
\label{constraintincomplete}
\end{align}

where $F \in \mathcal{H}^{\vee}$ such that $F(a + \langle b, \cdot \rangle) = a + \langle b, x_*\rangle$. The more general loss function we consider now is:
\begin{align}
S[\phi, \rho_2] &= -\int_{\mathbb{R}^n} \rho_2(y) \log\K{(}{)}{\int_{\mathbb{R}^n}\exp\K{(}{)}{\phi(x) + \langle y - x, \nabla \phi(x)\rangle}d\mu_1(x)}d^n y  + \int_{\mathbb{R}^n} \phi(x)d\mu_1(x)\\
&+ \int_{\mathbb{R}^n} \rho_2(x) \log\K{(}{)}{\frac{\rho_2(x)}{\rho_2^{ref}(x)}}d^n x,
\end{align}

where $\mu_2$ satisfies (\ref{constraintincomplete}) and $\rho_2^{ref} \in L^1(\mathbb{R}^n)$ is a reference probability density of choice with support $\mathrm{supp}(\rho)_2^{ref} = \Omega$ open and convex. We will denote the corresponding probability measure by $\mu_2^{ref}$. We pass to the dual:
\begin{align}
H[h, \phi] := \sup_{\rho_2, \int_{\mathbb{R}^n} \rho_2(y)d^n y = 1} \K{[}{]}{-\int_{\mathbb{R}^n} h(y)\rho_2(y)d^n y + F[h]- S[\phi, \rho_2] }.
\end{align}

The supremum is attained at:
\begin{align}
\rho_2(y) = \frac{1}{Z[h, \phi]}\rho_2^{ref}(y)\exp\K{(}{)}{-h(y)}\int_{\mathbb{R}^n}\exp\K{(}{)}{\phi(x) + \langle y - x, \nabla \phi(x)\rangle}d\mu_1(x),
\end{align}

where:
\begin{align}
Z[h, \phi] = \int_{\mathbb{R}^n\times \mathbb{R}^n } \rho_2^{ref}(y) \exp\K{(}{)}{-h(y) + \phi(x) + \langle y - x, \nabla \phi(x)\rangle}d\mu_1(x)d^n y.
\end{align}

In particular this means that if a minimimum to $\mathcal{H}$ exists, then the resulting $d\nu$ must be of the form
\begin{align}
d\nu(y) := \rho_2^{ref}(y)\exp\K{(}{)}{-h(y)}d^ny, \s h \in \mathcal{H}.
\end{align}

Substituting for the optimal $\rho_2$ we obtain
\begin{align}
H[h, \phi] = \log\K{(}{)}{Z[h, \phi]} + F[h] - \int_{\mathbb{R}^n} \phi(x)d\mu_1(x).
\end{align}

Define $\phi_h$ as the partial minimizer of $H$ over $\phi$ for fixed $h$, then $\phi_h$ is the Legendre transform of
\begin{align}
\psi_h(x) := \log\K{(}{)}{\int_{\mathbb{R}^n}\rho_2^{ref}(y)\exp\K{(}{)}{-h(y) + \langle k, y\rangle}d^n y},
\end{align}

and\footnote{It follows in particular:
\begin{align}
Z(h, \phi_h) = 1.
\end{align}
}

\begin{align}
G[h] := H[h, \phi_h] =  F[h] - \int_{\mathbb{R}^n} \phi_h(x)d\mu_1(x),
\end{align}

which is convex in $h$. Moreover we have
\begin{align}
\label{gradientphi}
\delta_{h(y)}\phi_h(x) &= -\delta_{h(y)}\K{.}{|}{\psi_h(k)}_{k = \nabla \phi_h(x) }\\
&= -\delta_{h(y)}\K{.}{|}{\log\K{(}{)}{\int_{\mathbb{R}^n}\rho_2^{ref}(y)\exp\K{(}{)}{-h(y) + \langle k, y\rangle}d^n y}}_{k = \nabla \phi_h(x)}\\
&=\rho_2^{ref}(y)\exp\K{(}{)}{-h(y) +\phi_h(x) + \langle y - x, \nabla \phi_h(x)\rangle}\\
&=K_h(y|x),
\end{align}

where by $K_h(y|x)$ we have denoted the Hessian martingale Kernel associated to $h$. Hence the local minimum of $G$, if it exists, would be attained where
\begin{align}
F[h] = \int_{\mathbb{R}^n}K_h(y|x)d\mu_1(x).
\end{align}

At this point we come to the main result of this section. In particular we will establish sufficient conditions on $\mathcal{H}$ such that $G$ is coercive and hence has a local minimum. 

\textbf{Theorem 2:} Let $\mathcal{H} \subset \mathcal{E}(\mathbb{R}^n)$ be a finite dimensional vector space closed under lower convex envelopes relative to $\Omega:= \mathrm{supp}(\rho_2^{ref})$ open and convex. Moreover let $\mathcal{H}$ contain all affine functions as well as functions whose lower convex envelope is strictly convex. Let $F \in \mathcal{H}^*$ and $\mu_1$ be a Radon probability measure with support in $\Omega$, such that
\begin{align}
&\int_{\mathbb{R}^n} \phi_h(x)d\mu_1(x) < \infty, \s \forall h \in (\mathcal{H} \cap \mathcal{K} / \mathcal{A}_{\Omega}) \backslash \{0\},\\
&F[h] > \int_{\mathbb{R}^n}h(x)d\mu_1(x) \s \forall h \in (\mathcal{H} \cap \mathcal{K} / \mathcal{A}_{\Omega}) \backslash \{0\},\\
&F[h] > 0 \s \forall h \in \mathcal{H}\backslash\{0\}\; with \; h(x) \geq 0 \s \forall x \in \mathbb{R}^n,
\end{align}

where $\mathcal{K} \subset \mathcal{E}(\mathbb{R}^n)$ is the cone of convex functions in $\mathcal{E}(\mathbb{R}^n)$, and $\mathcal{A}_{\Omega} \subset \mathcal{E}(\mathbb{R}^n)$ is the subspace of functions affine when restircted to $\Omega$. 
Then, the convex function
\begin{align}
G : \mathcal{H} / \mathcal{A}_{\Omega} &\rightarrow \mathbb{R},\\
h &\mapsto F[h] - \int_{\mathbb{R}^n} \phi_h(x)d\mu_1(x),
\end{align}

is coercive and hence has a local minimum.

\textbf{Proof:} First of all it is clear that $G$ is well defined on the quotient $\mathcal{H}/\mathcal{A}_{\Omega}$. Then for any $h \in (\mathcal{H}/\mathcal{A}_{\Omega}) \backslash \{0\}$ we have by Proposition 2,
\begin{align}
\lim_{r \rightarrow \infty}\frac{1}{r}G[rh] &= F[h] - \int_{\mathbb{R}^n} \lim_{r \rightarrow \infty} \frac{1}{r}\phi_{rh}(x)d\mu_1(x)\\
&= F[h] - \int_{\mathbb{R}^n} \mathrm{conv}_{\Omega}(h)(x) d\mu_1(x)\\
\label{convtrick}
&= (F[h - \mathrm{conv}_{\Omega}(h)]) + \K{(}{)}{F[\mathrm{conv}_{\Omega}(h)] - \int_{\mathbb{R}^n} \mathrm{conv}_{\Omega}(h)(x) d\mu_1(x)},
\end{align}

where $\mathrm{conv}_{\Omega}(h) \in \mathcal{E}(\Omega)$ and hence $\mathrm{conv}_{\Omega}(h) \in L^1_{\mu_1}(\mathbb{R}^n)$. Moreover we used the closedness of $\mathcal{H}$ under lower convex envelopes in order to evaluate $F$ in (\ref{convtrick}). Now, if $h$ is strictly convex, then $h = \mathrm{conv}_{\Omega}(h)$ over $\Omega$, hence the first term on the right hand side vanishes, however the second term is strictly greater than zero. Instead, if $h$ is not convex over $\Omega$, then $h - \mathrm{conv}_{\Omega}(h) \neq 0$ and $h(x) - \mathrm{conv}_{\Omega}(h)(x) \geq 0 \s \forall x \in \Omega$, hence the first term is strictly greater than zero. Therefore:
\begin{align}
\lim_{r \rightarrow \infty}\frac{1}{r}G[rh] > 0.
\end{align}

That is $G$ is coercive and hence has a local minimum.
\begin{align*}
\s\s\s\s\s\s\s\s\s\s\s\s\s\s\s\s\s\s\s\s\s\s\s\s\s\s\s\s\s\s\s\s\s\s\s\s\s\s\s\s\s\s\s\s\s\s\s\s\s\s\s\s\s\s\s\s\s\s\s\s\s\s\s\s\s\s\s \square
\end{align*}

We obtain the following corollary from the proof of Theorem 2, which can also be viewed as a generalization of the latter.

\textbf{Corollary:} If the vector space $\mathcal{H}$ is not closed under lower convex envelopes, then a sufficient condition for $G$ to have a local minimium is
\begin{align}
\sup_{f \in \mathcal{K}_h } \K{(}{)}{ F_2[h] - F_1 [f] } > 0 \s \forall h \in \mathcal{H},
\end{align}

where:
\begin{align}
\mathcal{K}_h := \{f \in \mathcal{K} \cap \mathcal{H} \, | \, f \geq \mathrm{conv}_{\Omega}(h)\}.
\end{align}

While such a condition only depends on the available data, in practice it is hard to verify.

\subsection{Global Calibration with Incomplete Marginals}\label{globalcalibration}

In the previous section we analyzed the variational problem of minimizing mutual information with linear constraints thereby obtaining $K$ and thus $\mu_2$ given $\mu_1$. We showed that a solution to such a variational problem exists for all $\mu_1$ satisfying partial no-arbitrage conditions, provided the vector space of calibration targets $\mathcal{H}$ has an additional property. Namely that of being closed under lower convex envelopes. While this property is enjoyed if $\mathcal{H} = \mathcal{E}(\mathbb{R}^n)$, in practice $\mathcal{H}$ will fail to satisfy it. Instead for arbitrary finite dimensional $\mathcal{H}_k \subset \mathcal{E}(\mathbb{R}^n)$, $k=1, \dots, N$, we will now show, under suitable assumptions on $\mu_0$, that the simultaneous variational problem for $\mu_1, \dots, \mu_N$ together with the Hessian Martingale Kernels between them, has a unique solution where here $\mu_i$ stands for the marginal at time $t_i$. We start from the full action:
\begin{align}
S[\phi_{\bullet}, \rho_{\bullet}] = \sum_{k = 0}^{N-1}&\K{[}{.}{-\int_{\mathbb{R}^n} \rho_{k+1}(y) \log\K{(}{)}{\int_{\mathbb{R}^n}\exp\K{(}{)}{\phi_{k+1}(x) + \langle y - x, \nabla \phi_{k+1}(x)\rangle}d\mu_k(x)}d^n y  }\\
&\K{.}{]}{+ \int_{\mathbb{R}^n} \phi_{k+1}(x)d\mu_k(x) + \int_{\mathbb{R}^n} \rho_{k+1}(x) \log\K{(}{)}{\frac{\rho_{k+1}(x)}{\rho_{k+1}^{ref}(x)}}d^n x},
\end{align}

where  $\rho_{k+1}^{ref} \in L^1(\mathbb{R}^n)$ is a reference probability density of choice with support $\mathrm{supp}(\rho)_{k+1}^{ref} = \Omega_{k+1}$ open and convex. Notice how now the previous densities $\rho_k$ are optimized for. $S$ is convex in $\rho_{\bullet}$ and concave in $\phi_{\bullet}$. We define the generalization
\begin{align}
H[\phi_{\bullet}, h_{\bullet}] = \sup_{\rho_{\bullet}, \int_{\mathbb{R}^n} \rho_{\bullet}(y)d^n y = 1} \K{[}{]}{\sum_{k = 0}^{N-1}\K{[}{]}{-\int_{\mathbb{R}^n} h_{k+1}(y)\rho_{k+1}(y)d^n y + F_{k+1}(h_{k+1})}- S[\phi_{\bullet}, \rho_{\bullet}] }.
\end{align}

Solving for $\rho_{k}$, $k = 1, \dots, N$ we obtain
\begin{align}
\rho_{k}(y) = \frac{1}{Z_k[h_{k},..,h_N, \phi_{k}, \dots, \phi_N]}&\rho_{k}^{ref}(y)\exp\K{(}{)}{-h_{k}(y) - \phi_{k+1}(y) + \int_{\mathbb{R}^n} d^n u K_{k+1}(u|y)}\\
&\cdot\int_{\mathbb{R}^n}\exp\K{(}{)}{\phi_{k}(x) + \langle y - x, \nabla \phi_{k}(x)\rangle}d\mu_{k-1}(x).
\end{align}

With the boundary conditions:
\begin{align}
&\phi_{N + 1} = 0,\\
&d\mu_0(x) \s \mathrm{given}.
\end{align}

We can now solve for the infimum over $\phi_{\bullet}$ iteratively starting from $\phi_N$ with $\phi_{N-1}, \dots, \phi_{1}$ fixed to obtain the following recursive identities. Denote
\begin{align}
&d\nu_{k}(y) = \rho_k^{ref}(y)\exp\K{(}{)}{-f_k(y)}d^n y,\\
&\phi_{k+1, f_{k+1}}(y) := \phi_{k+1}(y).
\end{align}

Then
\begin{align}
&f_k(y) = h_k(y) + \phi_{k+1, f_{k+1}}(y),\\
&f_{N+1} = 0.
\label{recursiveidentities}
\end{align}

Finally
\begin{align}
G[h_{\bullet}] &:= \inf_{\phi_{\bullet}} H[h_{\bullet}, \phi_{\bullet}]\\
&= \sum_{k = 1}^N F_{k}[h_k] -  \int_{\mathbb{R}^n}\phi_{1, f_1}(x)d\mu_0(x).
\end{align}
Following the recurisive identities (\ref{recursiveidentities}) and (\ref{gradientphi}), the local and global minimum of $G$, if it exists, is then attained where:

\begin{align}
F_k[h_k] &= \int_{\mathbb{R}^{n(k+1)}}h_k(x_k)\prod_{l=1}^N \K{(}{)}{K_{l, f_l}(x_l|x_{l-1})d^n x_l }d\mu_0(x_0)\\
&=\int_{\mathbb{R}^{n}}h_k(x_k)\rho_k(x)d^n x, \s \forall h_k \in \mathcal{H}_k, \s \forall k = 1, \dots, N.
\end{align}

\textbf{Theorem 3: } Consider $\bigoplus_{k=1}^N \mathcal{H}_k, \s \mathcal{H}_k \subset \mathcal{E}(\mathbb{R}^n) \s \forall k = 1, \dots, N$  with each $\mathcal{H}_k$ finite dimensional and containing all affine functions in $\Omega_k$. Let moreover $\mu_0$ have compact support inside $\Omega_1$\footnote{This includes the standard case of interest of a Dirac measure centered at a point $x_* \in \Omega_1$. We also remark that the compactness of the support of $\mu_0$ is not strictly necessary, in fact it suffices if $\mu_0$ is such that $\phi_{1, f_1} \in L_{\mu_0}(\mathbb{R}^n)$ for all choices of $h_k$ such that $\phi_{1, f_1}$ is finite. For $\Omega_0 = \mathbb{R}^n$, the latter condition is equivalent to $\mu_0 \in \mathcal{E}(\mathbb{R}^n)^*$.}. Then the convex function
\begin{align}
G_N: \bigoplus_{k=1}^N \mathcal{H}_k / \mathcal{A}_{\Omega_{k}} & \rightarrow \mathbb{R},\\
h_{\bullet} &\mapsto \sum_{k = 1}^N F_{k}[h_k] - \int_{\mathbb{R}^n}\phi_{1, f_1}(x)d\mu_0(x),
\end{align}

is coercive and hence has a global minimum if and only if  there are $\mu_0 \preceq \mu^R_1 \dots \preceq \mu^R_N$ Radon probability measures in the dual to $\mathcal{E}(\mathbb{R}^n)$ such that for all $k = 0, \dots, N - 1$
\begin{align}
&F_{k+1}(h) = \int_{\mathbb{R}^n}h(y)d\mu^R_{k+1}(y),\\
\label{ineqglobalconvex}
& \int_{\mathbb{R}^n}f(y)(d\mu^R_{k+1}(y) - d\mu^R_{k}(y)) > 0 \s \forall f \in (\mathcal{K} / \mathcal{A}_{\Omega_{k+1}}) \backslash \{0\}.\\
&  \mu^R_k \sim \mu_k^{ref} 
\end{align}

\textbf{Proof:}
\textbf{If:}
We first remark the straightforward fact that the existence of $\mu_0 \preceq \mu^R_1 \dots \preceq \mu^R_N$ Radon probability measures with $\mathrm{supp}(\mu^R_k) = \Omega_k$ (open and convex) implies that
\begin{align}
\Omega_0 \subset \Omega_1 \subset \dots \subset \Omega_N.
\end{align}

In particular, from Proposition 2 it follows that for $f \in \mathcal{E}$, $\mathrm{conv}_{\Omega_k}(f) \in \mathcal{E}(\Omega_l) \subset L^1_{\mu_l^R}(\mathbb{R}^n)$ for all $0 \leq l \leq k$, or it is identically $-\infty$. 
Thanks to the compact support of $d\mu_0$, by applying Proposition 2 iteratively, it follows
\begin{align}
\label{limitcoerciveglobal}
\lim_{r \rightarrow \infty}\int_{\mathbb{R}^n}\frac{1}{r} \phi_{1, r f_1}^*(x)d\mu_0(x) &= \int_{\mathbb{R}^n}\mathrm{conv}_{\Omega_1}(g_1^{\infty})(x)d\mu_0(x),
\end{align}

where $g_{n}^{\infty}$ is defined recursively as
\begin{align}
\label{ginftyrecursion1}
g_k^{\infty}(x) &= h_k(x) + \mathrm{conv}_{\Omega_{k+1}}(g_{k+1}^{\infty})(x),\\
\label{ginftyrecursion2}
g_{N+1}^{\infty}(x) &= 0.
\end{align}

Since $\mu_{k-1}^R \preceq \mu_{k}^{R}$ with $\mathrm{supp}(\mu_k) = \Omega_k$ for all $k = 1, \dots, N$ satisfying the strict inequalities (\ref{ineqglobalconvex}), it follows that for all $f \in (\mathcal{E}(\mathbb{R}^n) / \mathcal{A}_{\Omega_k}) \backslash{\{0\}}$
\begin{align}
\int_{\mathbb{R}^n}f(x)d\mu_k^R(x) > \int_{\mathbb{R}^n} \mathrm{conv}_{\Omega_k}(f)(x)d\mu_{k-1}^R(x).
\end{align}
Now, choose $h_{\bullet} \neq 0$, then there is $1\leq l \leq N$ such that $h_l$ is not affine in $\Omega_l$. Hence by induction,  (\ref{limitcoerciveglobal}), (\ref{ginftyrecursion1}),  (\ref{ginftyrecursion2})  lead to
\begin{align}
\lim_{r \rightarrow \infty} \frac{1}{r}G[rh_{\bullet}] &= \lim_{r \rightarrow \infty} \frac{1}{r}\K{(}{)}{ \sum_{k = 1}^N F_k(rh_k) - \int_{\mathbb{R}^n}\phi_{1, rf_1}(x)d\mu_0(x) }\\
&= \sum_{k = 1}^N F_k[h_k] - \int_{\mathbb{R}^n}\mathrm{conv}_{\Omega_1}(g^{\infty}_1)(x)d\mu_0(x)\\
&> \sum_{k = l+1}^N F_k[h_k] - \int_{\mathbb{R}^n}\mathrm{conv}_{\Omega_{l+1}}(g^{\infty}_{l+1})(x)d\mu_l^R(x)\\
&\geq 0.
\end{align}

\textbf{Only If:} Follows from the fact that the resulting optimal kernels $k-1 \rightarrow k$ are mutually absolutely continuous to $\mu_k^{ref}$.
\begin{align*}
\s\s\s\s\s\s\s\s\s\s\s\s\s\s\s\s\s\s\s\s\s\s\s\s\s\s\s\s\s\s\s\s\s\s\s\s\s\s\s\s\s\s\s\s\s\s\s\s\s\s\s\s\s\s\s\s\s\s\s\s\s\s\s\s\s\s\s \square
\end{align*}

\textbf{Remark:} While it was implicit in the proof of Theorem 3, it is instructive to observe explicitly from (\ref{limitcoerciveglobal}) an example of how the coercivity of $G$ requires the inequalities (\ref{ineqglobalconvex}). Consider the case where we have $h \in \mathcal{H}_k \cap \mathcal{H}_{l}$ convex and $k < l$. Choosing $h_k=-h, h_l=h$ with all others vanishing, equation (\ref{limitcoerciveglobal}) reduces to:
\begin{align}
\lim_{r \rightarrow \infty}\int_{\mathbb{R}^n}\frac{1}{r} \phi_{1, r f_1}^*(x)d\mu_0(x) &= \int_{\mathbb{R}^n}\mathrm{conv}_{\Omega_k}(-h + \mathrm{conv}_{\Omega_{l}}(h))d\mu_0(x)\\
&=\int_{\mathbb{R}^n}\mathrm{conv}_{\Omega_k}(-h + h)d\mu_0(x) = 0.
\end{align}

It follows
\begin{align}
\lim_{r \rightarrow \infty} \frac{1}{r}G[rh_{\bullet}] = -F_k[h] + F_l[h].
\end{align}
In other words $G$ coercive implies $F_k[h] < F_l[h]$.

\subsection{FX Market with Three Currencies}\label{fxexample}

Given three currencies $domestic$, $foreign_1$ and $foreign_2$, natural vanilla products used as calibration targets for the Risk-Neutral process on the exchange rates are FX options on each of the three exchange rates. Assuming no arbitrage and trivial interest rates for ease of exposition\footnote{Having deterministic rates does not change the analysis, while when rates are stochastic one needs to work first with the T-Forward measure and then pass to Risk-Neutral when simulating in Monte-Carlo.}, these will have prices
\begin{align}
C_{d \leftarrow f_1}(K, T) &:= \mathbb{E}_d\K{[}{]}{(X_{d \leftarrow f_1}(T) - K)^+},\\
C_{d \leftarrow f_2}(K, T) &:= \mathbb{E}_d\K{[}{]}{(X_{d \leftarrow f_2}(T) - K)^+},\\
C_{f_1 \leftarrow f_2}(K, T) &:= \mathbb{E}_{f_1}\K{[}{]}{(X_{f_1 \leftarrow f_2}(T) - K)^+}\\
&= \frac{1}{X_{d \leftarrow f_1}(0)}\mathbb{E}_{d}\K{[}{]}{(X_{d \leftarrow f_2}(T) - K \, X_{d \leftarrow f_1}(T))^+}.
\end{align}

In practice, at a given maturity $T$, we observe a finite number of payoffs, with a set of strikes $K_{ij}$, $i = 1,2,3$, $j = 1, \dots, M$. In this case, $h$ is of the form
\begin{align}
h_{a}(y_1, y_2) = \sum_{j = 1}^M a_{1j} (y_1 - K_{1j})^+ + \sum_{j = 1}^M a_{2j} (y_2 - K_{2j})^+ + \sum_{j = 1}^M a_{3j} (y_1 - K_{3j} y_2)^+,
\end{align}

and $y \in \Omega = \mathbb{R}_+^2$. Natural candidates for the reference marginals $\rho_k^{ref}$ are the product distributions of the densities obtained from Dupire's method applied to $X_{d\leftarrow f_1}$ and $X_{d\leftarrow f_2}$ separately. Note that the thus defined $\mathcal{H}$ fails in general to be closed under lower convex envelopes. In practice lower convex envelopes are well approximated by convex functions in $\mathcal{H}$, thus failure to calibrate sequentially to market prices would be due to a very tight no-arbitrage scenario. On the other hand global calibration across all maturities is guaranteed by Theorem 3.

\section{Monte Carlo Calibration of Hessian Martingales}
In this section we will use the machinery previously developed in order to sketch a strategy to fit a multi-dimensional continuous martingale to prices of a set of European options corresponding to a grid of strikes and maturities. Examples of such payoffs are basket options of the form:
\begin{align}
f_{(a, b)}(x) := \max\K{(}{)}{0, a + \langle b, x\rangle}, \s a \in \mathbb{R}, b \in \mathbb{R}^n.
\end{align}

In particular we focus on a Monte-Carlo pricing framework. At $t = T_m$ we will have obtained an empirical distribution $\{x_1, \dots, x_N\} \in \mathbb{R}^n$. The aim now is to diffuse these to $t = T_{m+1}$ as follows\footnote{or, for manifestly positive domains
\begin{align}
\log(x_k^i) \mapsto \log(x_k^i) - \frac{1}{2(x_k^i)^2}(\Sigma_k \Sigma_k^T)^{ii} + \frac{1}{x_k^i}(\Sigma_k u_k)^i,
\end{align}

maybe subdividing it into $n$ steps in between by scaling $\Sigma$ with $(T_{m+1} - T_m)/n$ for each step.
}
\begin{align}
x_k \mapsto x_k + \Sigma_k u_k,
\end{align}

where $\Sigma_k$ is a volatility matrix at $x_k$, $u_k$ is a $n$-dimensional standard normal random variable. From the resulting distribution we obtain prices of the vanillas maturing at $T_{m+1}$. The aim is to fit $\Sigma_k$ to the market prices, which can then be achieved iteratively. The question here is how to model $\Sigma_k$. Given 
what we have developed so far, without loss of generality (assuming relatively small diffusion times $\Delta T := T_{m+1} - T_m$) we can choose $\Sigma_k$ such that
\begin{align}
\Sigma_k \Sigma_k^T &= g^{-1}(x_k) = \dfrac{\int_{\mathbb{R}^n} (y - x_k)(y - x_k)^T \exp\K{(}{)}{\langle \nabla \phi(x_k), y\rangle}d\nu(y)}{\int_{\mathbb{R}^n} \exp\K{(}{)}{\langle \nabla \phi(x_k), y\rangle}d\nu(y)}\\
&= H(\psi)(\nabla \phi(x_k)).
\end{align}

Therefore, given a choice of $d\nu(y)$ and given $x_k$, we obtain $k(x_k) := \nabla \phi(x_k)$ by determining the maximum over $k$ of:
\begin{align}
\langle x, k\rangle - \psi(k),
\end{align}

which is easily achieved by steepest descent. In principle we would have to solve $N$ such optimization problems (one for each simulated point $x_k$). However, in practice, we can do this for a different set of carefully chosen knots $u_1, \dots, u_{N_{knots}} \in \mathbb{R}^n$ with $N_{knots} << N$ and interpolate to obtain an estimate at $x_k$. More precisely, we will interpolate each component of the $\Sigma$ matrix. What is left to determine is $d\nu_{m}(y)$ at maturity $T_{m}$ which as shown in section \ref{IncompleteMarginals} is of the form given recursively backwards across maturities by:
\begin{align}
d\nu_{m}(y) = \rho_m^{ref}(y)\exp\K{(}{)}{-h_{m}(y) - \phi_{m+1}(y)}d^n y, \s h_{m} \in \mathcal{H}_m.
\end{align}

Here $\rho_m^{ref}(y)$ is a reference density as detailed in Theorem 3 (e.g., the product density resulting from Dupire's approach in each dimension separately). Instead $\mathcal{H}_m$ denotes the vector space of (discounted) European payoffs maturing at $T_m$, and $\phi_{m+1}$ results from the step at maturity $T_{m+1}$.
More explicitly we can choose a basis $h_{im}$ for such payoffs at $T_m$, then $h_m$ will be of the form:
\begin{align}
h_{m} = \sum_i a_{im} h_{im}.
\end{align}

The coefficients $a_{im}$ are optimized so as to hit market prices at all maturities $T_1, \dots, T_N$.

Once the process is calibrated, we will have obtained the $\Sigma(T_m)$ interpolants at every maturity $T_m$, and we can then simulate paths at intermediate time-points by using a Brownian bridge with the given $\Sigma$ between the simulated points at $T_m$ and $T_{m+1}$.
\newpage
\section{Conclusion}

In the present paper, we constructed what we called ``Hessian martingales", namely the unique ``most unbiased" (Markov) martingales generating marginals devoid of calendar arbitrage. Moreover we outlined a concrete method to construct a unique Hessian martingale from arbitrage free prices of a finite set of European contingent claims. The  theoretical construction is clearly outlined in the case of deterministic (discounting) interest rates, but it can in principle be extended to the case of stochastic interest rates, either by including bonds as tradeables or by imposing a precalibrated rates dynamics with some assumptions on the correlation structure between interest rates and the underlyings in question. The former approach is elegant and immediate from a theoretical standpoint, but impractical as it quickly turns into a very high dimensional problem. It is therefore interesting to explore the extension of the present construction to incorporate pre-calibrated stochastic exogenous dynamics, as e.g. the discounting rate. From a more theoretical/technical standpoint it would be fruitful to fully understand the construction in continuous time, which here was only briefly sketched. We leave to future work a study of various implementations and examples of the approach presented here. 

\medskip
\medskip
\medskip

\textit{Acknowledgements} I am grateful to Gianluca Lauteri for early stimulating discussions, and for helpful comments throughout. I am also indebted to Francesco Mattesini for a careful reading of an early version of the manuscript and in particular for pointing out the related recent work in \cite{Wiesel},  as well as to Alexander Igamberdiev and Walter Reusswig for useful comments that helped improve the clarity of exposition.

\newpage
\appendix
\section{Failure of local correlation models}

A popular class of models to fit prices of European basket options that depend on multiple assets (e.g. stocks, FXs...) was studied in \cite{Guyon13}. Such a class consists of continuous local Markov martingales with a \textsl{locally gaussian Markov copula}. For the case of $n$ assets, by this we mean processes with an associated forward Kolmogorov equation of the form
\begin{align}
\partial_t \rho_t(x) = \frac{1}{2}\partial_i \partial_j \K{(}{)}{\eta_t^{ij}(x)\rho_t(x)},
\end{align}

where:
\begin{align}
\eta_t^{ii} = \eta_t^{ii}(x^i).
\end{align}

In other words, these processes are continuous local martingales with the extra condition that for each time $t$, the diagonal elements of the covariance matrix $\eta_t^{-1}$ only depend on their corresponding underlying. The off-diagonal elements may instead depend on all the underlyings (however constrained to ensure the positive semi-definiteness of $\eta_t^{-1}$).
In this section we will show how this class of models is not sufficient to reproduce arbitrary arbitrage-free European basket option prices. In fact we prove the stronger result implying that, e.g. in the context of the FX market, the aforementioned class of models is insufficient in general to calibrate to arbitrage free prices of FX options relative to three currencies (see section \ref{fxexample}).

\textbf{Proposition 3:} Consider the functions in $\mathcal{E}(\mathbb{R}^2)$ given by
\begin{align}
f_K(x,y) := \max(x - Ky, 0), \s K \in \mathbb{R}_+.
\end{align}

then, there is a continuous Markov martingale with marginals $\rho_t(x,y)$ and covariance matrix $g_t^{-1}$ that does not admit any (associated) continuous Markov martingale with a locally Gaussian Markov copula and covariance matrix $\eta_t^{-1}$, whose marginals $\widetilde{\rho}_t$ have the same 1-d marginals as $\rho_t$ and the same expectation values of $f_K$ for all $K \in \mathbb{R}_+$ and $t$ in an open interval.

\textbf{Proof:}

Let $g_t^{-1}$ denote the covariance matrix at time $t$ of the true process and be such that $\rho_t$ has support in $\mathbb{R}_+^2$. Let
\begin{align}
C(t, K) &:= \int_{\mathbb{R}_+^2} f_K(x,y) \rho_t(x,y) dx dy.
\end{align}

Then
\begin{align}
2\partial_t C(t, K) &= \int_{\mathbb{R}_+} \K{[}{]}{g_t^{11}(Ky,y) - 2 K g_t^{12}(Ky,y) + K^2 g_t^{22}(Ky,y)}\rho_t(Ky,y) dy.
\end{align}

If an associated $\eta_t^{-1}$ exists with the requirements stated in the proposition, then it must equally hold
\begin{align}
2\partial_t C(t, K) &= \int_{\mathbb{R}_+} \K{[}{]}{\eta_t^{11}(Ky) - 2 K \eta_t^{12}(Ky,y) + K^2 \eta_t^{22}(y)}\widetilde{\rho}_t(Ky,y) dy.
\end{align}

It is straightforward to show that $\eta_t^{11}, \eta_t^{22}$ are fully determined by $\rho_t$ and by $g_t^{11}, g_t^{22}$ respectively:
\begin{align}
\eta_t^{11} = \sigma_{1, t}^2(x) &:= \int_{\mathbb{R}_+} g_t^{11}(x,y)\rho_{t}(y|x)dy,\\
\eta_t^{22} = \sigma_{2, t}^2(y) &:= \int_{\mathbb{R}_+} g_t^{22}(x,y)\rho_{t}(x|y)dx.
\end{align}

The positive semi-definiteness of $\eta_t^{-1}$, is equivalent to the following two set of inequalities:
\begin{align}
\label{ineq1}
2\partial_t C_t(t, K) &\geq \int_{\mathbb{R}_+} \K{[}{]}{\sigma_{1,t}(Ky) - K \sigma_{2,t}(y)}^2  \widetilde{\rho}_t(Ky, y) dy,\\
\label{ineq2}
2\partial_t C_t(t, K) &\leq \int_{\mathbb{R}_+} \K{[}{]}{\sigma_{1,t}(Ky) + K \sigma_{2,t}(y)}^2  \widetilde{\rho}_t(Ky, y) dy.
\end{align}

Now we proceed to construct a counterexample $g_t^{-1}$. We consider a process that starts at $(x,y)=(1,1)$ at $t=0$ and proceeds as two independent geometric Brownian motions until $t=1$. After $t = 1 + \epsilon$ for $\epsilon > 0$, instead, the off-diagonal elements remain null, while the diagonal elements crystallize to a new form as follows:
\begin{align}
g_{t, \epsilon}^{12}(x, y) &= g_{t, \epsilon}^{21}(x, y) = 0,\\
g_{t, \epsilon}^{11}(x, y) &= x^2\K{(}{)}{(1 - \lambda_{\epsilon}(t)) + \lambda_{\epsilon}(t) a y^{2\alpha}},\\
g_{t, \epsilon}^{22}(x, y) &= y^2.
\end{align}

where we have made the dependence of $g_{t}^{-1}$ on $\epsilon$ explicit as $g_{t, \epsilon}^{-1}$. Here $\alpha \in \mathbb{R}$, and $\lambda_{\epsilon}$ is a smooth, monotonically increasing mollified characteristic function of the interval $[1 + \epsilon, \infty)$. Moreover, $\lambda_{\epsilon}(t) = 1$ for $t \geq 1 + \epsilon$ and $\lambda_{\epsilon}(t) = 0$ for $t \in [0,1]$. We shall take the family $\lambda_{\epsilon}$ to be smooth in $\epsilon$ (e.g. by rescaling). The process defined by $g_t$ is a true martingale for all $\alpha \in \mathbb{R}$. We shall fix the constant $a$ such that:
\begin{align*}
\int_{\mathbb{R}_+}a y^{2\alpha} \rho_1(y|x)dy = 1,
\end{align*}

where $\rho_1$ is the density at $t=1$ resulting from $g_{t, \epsilon}$. By construction it is independent of $\epsilon$ and given by
\begin{align}
\rho_1(x,y) = \psi(x)\psi(y),
\end{align}

with:
\begin{align}
\psi(x) := \frac{1}{\sqrt{2\pi}}\frac{1}{x}\exp\K{(}{)}{-\frac{1}{2}\K{(}{)}{\log(x) + \frac{1}{2}}^2}.
\end{align}

Therefore
\begin{align*}
a = e^{-\alpha(2\alpha-1)}.
\end{align*}

It follows that
\begin{align}
\sigma^2_{1, 1+\epsilon, \epsilon}(x) &:= \eta_{1 + \epsilon, \epsilon}^{11} = x^2 + \delta(\epsilon, x),\\
\sigma^2_{2, 1 + \epsilon, \epsilon}(y) &:= \eta_{1 + \epsilon, \epsilon}^{22} = y^2,
\end{align}

where $\delta(\epsilon, x)$ vanishes for all $x$ as $\epsilon \rightarrow 0$. Let $C_{\epsilon}(t, K)$ denote the expectation values of $f_K$ induced by $g_{t, \epsilon}$. We now inspect (\ref{ineq1}) and (\ref{ineq2}) at $t = 1+\epsilon$ as $\epsilon \rightarrow 0$. The first becomes trivial, while the second becomes:
\begin{align}
2\partial_t C_0(1^+, K) \leq  4\int_{\mathbb{R}_+} K^2 y^2 \widetilde{\rho}_1(Ky, y)dy.
\end{align}

Let $f: \mathbb{R}_+ \rightarrow \mathbb{R}_+$, then:
\begin{align}
2\int_{\mathbb{R}_+} f(K)\partial_t C_0(1^+, K)dK  \stackrel{!}\leq 4 \int_{\mathbb{R}_+^2} f(K) K^2 y^2 \widetilde{\rho}_1(Ky, y)dy dK.
\end{align}

Changing variables, the above is equivalent to:
\begin{align}
\label{inequalityxy}
2\int_{\mathbb{R}_+} f\K{(}{)}{K}\partial_t C_0(1^+, K)dK  \stackrel{!}\leq 4 \int_{\mathbb{R}_+^2} f\K{(}{)}{\frac{x}{y}} x^2 \widetilde{\rho}_1(x, y)dx \frac{dy}{y}.
\end{align}

The difficulty we would like to overcome at this point lies in not knowing $\widetilde{\rho}_1$. However, if we can turn the r.h.s into an expectation value of a function solely of $x$ or of $y$, then by virtue of the fact that $\widetilde{\rho}_1$ has the same 1-d marginals as $\rho_1$, we can replace the former with the latter. The desired aim is achieved with either $f(K) = K^{-1}$ or $f(K) = K^{-2}$. Choosing the former, inequality (\ref{inequalityxy}) is satisfied, while choosing the latter, (\ref{inequalityxy}) becomes
\begin{align}
2\int_{\mathbb{R}_+} \partial_t C_0(1^+, K)\frac{dK}{K^2} &\stackrel{!}\leq 4 \int_{\mathbb{R}_+^2} y  \,\widetilde{\rho}_1(x, y)dx dy\\
&= 4 \int_{\mathbb{R}_+^2} y  \, \rho_1(x, y)dx dy\\
&= 4.
\end{align}

Instead, inserting $g_{1^+}^{-1}$ into the formula for $\partial_t C(1^+, K)$ we obtain
\begin{align}
2\int_{\mathbb{R}_+} \partial_t C_0(1^+, K)\frac{dK}{K^2} &= \int_{\mathbb{R}_+^2} \frac{y}{x^2}\K{(}{)}{ax^2y^{2\alpha} + x^2} \rho_1(x,y)dx dy\\
&= a \int_{\mathbb{R}_+} y^{1 + 2\alpha} \psi(y)dy + 1\\
&= e^{2\alpha} + 1.
\end{align}

Choosing $2\alpha > \log(3)$, we reach a contradiction. 
\begin{align*}
\s\s\s\s\s\s\s\s\s\s\s\s\s\s\s\s\s\s\s\s\s\s\s\s\s\s\s\s\s\s\s\s\s\s\s\s\s\s\s\s\s\s\s\s\s\s\s\s\s\s\s\s\s\s\s\s\s\s\s\s\s\s\s\s\s\s\s \square
\end{align*}

\newpage
\section{Convex Envelopes}\label{convexenvelopes}

Here we detail the proofs of the propositions in section \ref{IncompleteMarginals}, which constitute standard results in convex analysis.

\textbf{Proof of Proposition 1:} 
\begin{align}
\mathrm{conv}_{\Omega}(f)(x) &= \sup \{a + \langle b,x \rangle \, | \, a + \langle b,y \rangle \leq f(y) \s \forall y \in \Omega, a \in \mathbb{R}, b \in \mathbb{R}^n \}\\
&=\sup_{b \in \mathbb{R}^n} \K{(}{)}{\inf_{y \in \Omega} \K{(}{)}{f(y) - \langle b,y\rangle} + \langle b, x\rangle}\\
&=\sup_{b \in \mathbb{R}^n} \K{(}{)}{\langle b, x\rangle - \sup_{y \in \Omega} \K{(}{)}{\langle b,y\rangle - f(y)}}\\
&=\sup_{b \in \mathbb{R}^n} \K{(}{)}{\langle b, x\rangle - f_{\Omega}^*(b)}
\end{align}

\begin{align*}
\s\s\s\s\s\s\s\s\s\s\s\s\s\s\s\s\s\s\s\s\s\s\s\s\s\s\s\s\s\s\s\s\s\s\s\s\s\s\s\s\s\s\s\s\s\s\s\s\s\s\s\s\s\s\s\s\s\s\s\s\s\s\s\s\s\s\s \square
\end{align*}

\textbf{Proof of Proposition 2:} 
We split the proposition into three parts: 
\begin{itemize}
\item[\textbf{1.}] $\phi_h \in C(\langle \Omega \rangle) \cup \{-\infty\}$, 
\item[\textbf{2.}] $\lim_{r \rightarrow \infty} \frac{1}{r}\phi_{rh} = \mathrm{conv}_{\Omega}(h)$, $\mathrm{conv}_{\Omega}(h) \in C(\langle \Omega\rangle) \cup \{-\infty\}$, 
\item[\textbf{3.}] $h \in \mathcal{E}(\Omega) \Rightarrow \mathrm{conv}_{\Omega}(h) \in \mathcal{E}(\Omega) \cup \{-\infty\}$.
\end{itemize}

\subsubsection*{Part 1}
For $x \in \langle \Omega\rangle$, we can choose $y_0, \dots, y_{n} \in \Omega$ and probability vector $p_0, \dots, p_n$ such that
\begin{align}
x = \sum_{i=0}^{n} p_i y_i.
\end{align}

Then
\begin{align}
f_{h,x}(k) &:= -\log\K{(}{)}{\int_{\Omega} \rho(y)\exp\K{(}{)}{\langle k, y - x\rangle} - h(y)}\\
&= -\sum_{i=0}^n p_i \log\K{(}{)}{\int_{\Omega} \rho(y)\exp\K{(}{)}{\langle k, y - y_i\rangle} - h(y)}
\end{align}

Let $v \in S^{n-1}$ such that $k = \|k\|v$.  Since $\Omega$ is open, for each $y_i$, there is $\delta_i > 0$ such that the open set
\begin{align*}
U_i := \{y \in \Omega \, | \, \langle v, y - y_i\rangle \geq \delta_i\}
\end{align*}

is not empty. Given that $h \in C(\Omega)$\footnote{It suffices in fact that $h \in L^1_{loc}(\Omega)$ and locally bounded above.} and that $\mathrm{supp}(\rho) = \Omega$,
\begin{align}
c_i := \int_{U_i} \rho(y)\exp\K{(}{)}{-h(y)} > 0.
\end{align}

Then
\begin{align}
f_{h,x}(k) \leq -\K{(}{)}{\sum_{i=0}^n p_i \delta_i} \|k\| - \K{(}{)}{\sum_{i=0}^n p_i \log(c_i)}.
\end{align}

Hence $f_{h,x}(k)$ is and anti-coercive (and concave) and consequently:
\begin{align}
\phi_h(x) = \sup_{k \in \mathbb{R}^n} f_{h,x}(k) \in \mathbb{R} \cup \{-\infty\}.
\end{align}

Therefore, since $\phi_h$ is convex it is either continuous in $\langle \Omega\rangle$ or identically equal to $-\infty$.
\subsubsection*{Part 2}
We can express
\begin{align}
\frac{1}{r}\phi_{rh}(x) = \sup_{k \in \mathbb{R}^n}\K{(}{)}{\langle x, k\rangle - \psi_r(k)},
\end{align}

where
\begin{align}
\psi_r(k) :=  \frac{1}{r}\log\K{(}{)}{\int_{\mathbb{R}^n}\rho(y)\exp\K{(}{)}{r (\langle k, y\rangle - h(y))}d^n y}.
\end{align}

Since $h \in C(\Omega)$ and $\mathrm{supp}(\rho) = \Omega$,
\begin{align}
\psi_{\infty}(k) &:= \lim_{r \rightarrow \infty}\psi_r(k)\\
&= \lim_{r \rightarrow \infty}\log\K{(}{)}{\K{\|}{\|}{\exp\K{(}{)}{\langle k, \cdot\rangle - h}}_{\rho, r}}\\
&= \sup_{y \in \Omega}\K{(}{)}{\langle k, y\rangle - h(y)}.
\end{align}

Let
\begin{align}
\Omega^{*} := \{k \in \mathbb{R}^n \, | \, |\psi_{\infty}(k)| < \infty\}.
\end{align}

Then $\psi_{\infty}$ is convex and therefore continuous over $\Omega^*$. For all compact sets $K \subset \Omega^*$, $\psi_r(k)$ converges uniformly to $\psi_{\infty}$. Namely for all $\epsilon > 0$, there is $R > 0$ such that for all $r > R$
\begin{align}
|\psi_{\infty}(k) - \psi_r(k)| < \epsilon \s \forall k \in K.
\end{align}

This follows from the fact that $\psi_r(k)$ is monotonically increasing and continuous in $r$ for $r > 0$. Monotonicty of $\psi_r$ follows from
\begin{align}
\partial_r \psi_r(k) = \frac{1}{r^2}\int_{\mathbb{R}^n} \sigma_{k,r}(y)\log\K{(}{)}{\frac{\sigma_{k,r}(y)}{\rho(y)}}d^n y \geq 0,
\end{align}

where
\begin{align}
\sigma_{k,r}(y) := \frac{\rho(y)\exp\K{(}{)}{r (\langle k, y\rangle - h(y))}}{\int_{\mathbb{R}^n}\rho(u)\exp\K{(}{)}{r (\langle k, u\rangle - h(u))}d^n u}.
\end{align}

We have:
\begin{align}
f_x(k) := \langle x,k \rangle - \psi_{\infty}(k) \leq \langle k, x - y\rangle + h(y) \s \forall y \in \Omega.
\end{align}

Now, let $v \in S^{n-1}$ such that $k = \|k\|v$, then for all $x \in \langle \Omega\rangle$ we can choose $\delta > 0$, $y_0, \dots, y_n \in \Omega$, and probability vector $p_0, \dots, p_n$ such that:
\begin{align}
x - \delta v = \sum_{i=0}^n p_i y_i.
\end{align}

It follows
\begin{align}
f_x(k) \leq - \|k\| \delta + c \sum_{i=0}^n p_i h(y_i).
\end{align}

That is, $f_x$ is concave and anti-coercive on $\mathbb{R}^n$, hence it is either equal to $-\infty$ for all $x \in \langle \Omega\rangle$, corresponding to the case $\Omega^* = \emptyset$, or it attains a maximum for all $x \in \langle \Omega \rangle$. In the latter case, corresponding to $\Omega^* \neq \emptyset$, it then follows from the uniform convergence of $\psi_r$ on compact sets $K \subset \Omega^*$ that for each $x \in \langle \Omega \rangle$, there exists a compact $K_x \subset \Omega^*$ such that
\begin{align}
\lim_{r \rightarrow \infty} \frac{1}{r}\phi_{rh}(x)  &= \lim_{r \rightarrow \infty}\sup_{k \in \mathbb{R}^n} (\langle x,k \rangle - \psi_{r}(k))\\
&= \lim_{r \rightarrow \infty}\sup_{k \in K_x} (\langle x,k \rangle - \psi_{r}(k))\\
&= \sup_{k \in K_x} (\langle x,k \rangle - \psi_{\infty}(k))\\
&= \sup_{k \in \mathbb{R}^n} (\langle x,k \rangle - \psi_{\infty}(k)).
\end{align}

Finally
\begin{align}
\label{convomega}
\lim_{r \rightarrow \infty} \frac{1}{r}\phi_{rh}(x) = \mathrm{conv}_{\Omega}(h)(x) \s \forall x \in \langle \Omega \rangle.
\end{align}

Since $\mathrm{conv}_{\Omega}(h)$ is finite and convex in $\langle \Omega \rangle$, it is continuous there. 
\subsubsection*{Part 3}
If $\Omega^* \neq \emptyset$, there exists a $k_0 \in \Omega^*$, therefore for all $x \in \mathbb{R}^n$
\begin{align}
\mathrm{conv}_{\Omega}(h)(x) \geq \langle k_0, x\rangle - \psi_{\infty}(k_0) \geq -\max(\|k_0\|, |\psi_{\infty}(k_0)|)(1 + \|x\|).
\end{align}

Moreover, for all $x \in \Omega$ it holds
\begin{align}
\mathrm{conv}_{\Omega}(h)(x) \leq h(x) \leq c(1 + \|x\|).
\end{align}

\begin{align*}
\s\s\s\s\s\s\s\s\s\s\s\s\s\s\s\s\s\s\s\s\s\s\s\s\s\s\s\s\s\s\s\s\s\s\s\s\s\s\s\s\s\s\s\s\s\s\s\s\s\s\s\s\s\s\s\s\s\s\s\s\s\s\s\s\s\s\s \square
\end{align*}

\newpage
\section{Why Markov?}

In this section we give a sketch of a proof that the (Markovian) Hessian martingales are the most unbiased given marginals, across all, not necessarily Markov, martingales\footnote{The argument generalizes in a straightforward way also to the case of incomplete marginals in the sense of section \ref{globalcalibration}.}. We have:

\textbf{Proposition 4:} Let $\mu_{i}$ with $i = 0, \dots N$ be Radon probability measures over $\mathbb{R}^n$ with $\mu_{i} \preceq \mu_{j}$, $i < j$. Consider the following functional on the space of joint Radon probability measures

\begin{align}
&S[\mu, b, c] := \int_{\mathbb{R}^{n(N+1)}} d\mu(x_0, ..., x_N) \log\K{(}{)}{\frac{d\mu(x_0, \dots, x_N)}{\prod_{k=0}^N d\mu_k(x_k)}} \\
&- \sum_{k = 0}^{N-1} \int_{\mathbb{R}^{n(N+1)}}\K{\langle}{\rangle}{ b_k(x_0, \dots, x_k), x_{k+1} - x_{k}} d\mu(x_0, \dots, x_N) \\
&- \sum_{k = 0}^{N} \int_{\mathbb{R}^{n(N+1)}} c_{k}(x_{k})\K{(}{)}{d\mu(x_0, \dots, x_N) - d\mu_k(x_k)},
\end{align}
where $b_k \in C^b(\mathbb{R}^{n(N+1)}, \mathbb{R}^n)$, $c_k \in \mathcal{E}(\mathbb{R}^n)$ are Lagrange multiplier functions imposing that $\mu$ is a martingale. It follows that if $\mu_0, \dots, \mu_N$ admit a Hessian martingale, then the resulting joint distribution is the minimum of  $S$.

\textbf{Proof Sketch:}

$S$ is strictly convex in $\mu$, hence if it has a local minimum, it will be of the form (we have absorbed an overall constant in the $c$ functions)

\begin{align}
d\mu(x_0, \dots, x_N) = \prod_{k = 0}^N d\mu_i(x_i) \exp\K{(}{)}{\sum_{k = 0}^{N} c_k(x_k) + \sum_{k = 0}^{N-1}\langle b_k(x_0, \dots, x_k), x_{k+1} - x_k\rangle}.
\end{align}
On the other hand, given optimal kernels, the corresponding $\mu$ is of the form

\begin{align}
d\mu(x_0, \dots, x_N)  =  \prod_{k = 0}^N d\mu_i(x_i)\exp\K{(}{)}{\sum_{k = 0}^{N-1} \K{[}{]}{\widetilde{c}_{k+1}(x_{k+1}) + \widetilde{a}_k(x_k) + \langle \widetilde{b}_k(x_k), x_{k+1}\rangle } },
\end{align}
which is of the optimal form for $S$ with $b_k = \widetilde{b}_k$ and

\begin{align}
c_0(x_0)  &= \widetilde{a}_0(x_0) + \langle \widetilde{b}_0(x_0), x_0\rangle = \phi_0(x_0)\\
c_{k+1}(x_{k+1}) &= \widetilde{c}_{k+1}(x_{k+1}) + \widetilde{a}_{k+1}(x_{k+1}) + \langle \widetilde{b}_{k+1}(x_{k+1}), x_{k+1}\rangle = \widetilde{c}_{k+1}(x_{k+1}) + \phi_{k+1}(x_{k+1}).
\end{align}

\begin{align*}
\s\s\s\s\s\s\s\s\s\s\s\s\s\s\s\s\s\s\s\s\s\s\s\s\s\s\s\s\s\s\s\s\s\s\s\s\s\s\s\s\s\s\s\s\s\s\s\s\s\s\s\s\s\s\s\s\s\s\s\s\s\s\s\s\s\s\s \square
\end{align*}

\section{Examples of Hessian Martingales}

In this section we revisit some classic examples of continuous martingales in the light of optimal kernels as well as construct some less well known ones. We start with revisiting well-known 1 dimensional examples:

\subsection{Dimension 1}

Here we revisit well known examples of 1d volatilities with skew $\sigma(x)$ and construct the corresponding finite time $\delta t$ kernels by vieweing the variance $\sigma(x) \sqrt{\delta t}$ as

\begin{align}
\partial_x^2 \phi(x) = \frac{1}{\sigma^2(x) \delta t}.
\end{align}
Then we proceed ``backwards" and reconstruct $\psi$, hence $d\nu$ and finally $d\mu(\, \cdot \, | \, \cdot \,)$ from $\sigma$.

\subsubsection{Normal}

Here $\Omega = \mathbb{R}$ and
\begin{align}
\sigma : \Omega &\rightarrow \mathbb{R}_{\geq 0},\\
x &\mapsto 1.
\end{align}

Therefore
\begin{align}
\phi(x) = \frac{1}{2 \delta t} x^2 + ax + b.
\end{align}
We choose $x_* = 0$, that is we fix $a,b$ such that $\phi(0) = \partial_x \phi(0) = 0$, hence $a = b = 0$. Then
\begin{align}
\psi(k) = \sup_{x \in \Omega} \K{(}{)}{kx - \phi(x)}.
\end{align}

The maximum is attained where
\begin{align}
k = \partial_x \phi(x) = \frac{x}{\delta t}.
\end{align}
Solving for $x$,

\begin{align}
\psi(k) = \frac{\delta t}{2} k^2.
\end{align}

Hence
\begin{align}
\int_{\mathbb{R}} d\nu(y) \exp\K{(}{)}{ky} &= \exp\K{(}{)}{\psi(k)}\\
&= \exp\K{(}{)}{\frac{\delta t}{2} k^2}.
\end{align}

In order to extract $d\nu$ we can proceed in several ways, but we opt for analytic continuation:
\begin{align}
\chi_{\nu}(k) &:= \int_{\mathbb{R}} d\nu(y) \exp\K{(}{)}{iky}\\
&= \exp\K{(}{)}{\psi(ik)}\\
&= \exp\K{(}{)}{-\frac{\delta t}{2} k^2}.
\end{align}
Therefore

\begin{align}
d\nu(y) &= \frac{dy}{2\pi} \int_{\mathbb{R}} \chi_{\nu}(k)\exp\K{(}{)}{-iky}d^n k\\
&= \frac{dy}{2\pi} \int_{\mathbb{R}}\exp\K{(}{)}{-\frac{\delta t}{2} k^2}\exp\K{(}{)}{-iky}d^n k\\
&= \frac{1}{\sqrt{2\pi \delta t}}\exp\K{(}{)}{-\frac{y^2}{2\delta t}}dy.
\end{align}

Finally
\begin{align}
d\mu(y|x) &= d\nu(y)\exp\K{(}{)}{\phi(x) + \partial_x \phi(x) (y - x)}\\
&= \frac{1}{\sqrt{2\pi \delta t}}\exp\K{(}{)}{-\frac{(y - x)^2}{2\delta t}}dy,
\end{align}

which is indeed what we expect.
\subsubsection{LogNormal}

Here $\Omega = (0, \infty)$ and
\begin{align}
\sigma : \Omega &\rightarrow \mathbb{R}_{\geq 0},\\
x &\mapsto x.
\end{align}

Therefore
\begin{align}
\phi(x) = -\frac{1}{\delta t} \log(x) + ax + b.
\end{align}

We choose $x_* = 1$, that is we fix $a,b$ such that $\phi(1) = \partial_x \phi(1) = 0$, hence
\begin{align}
\phi(x) = \frac{1}{\delta t} \K{(}{)}{-\log(x) + x - 1}.
\end{align}

Then
\begin{align}
\psi(k) = \sup_{x \in \Omega} \K{(}{)}{ kx - \phi(x)}.
\end{align}

The maximum is attained where
\begin{align}
k = \partial_x \phi(x) = \frac{1}{\delta t}\K{(}{)}{-\frac{1}{x} + 1}.
\end{align}
In particular $\Omega^* = \K{(}{)}{-\infty, \frac{1}{\delta t}}$. Solving for $x$:

\begin{align}
x = \frac{1}{1 - k \delta t}.
\end{align}

Therefore
\begin{align}
\psi(k) = - \frac{1}{\delta t} \log\K{(}{)}{1 - k \delta t}.
\end{align}

Hence
\begin{align}
\int_{\mathbb{R}} d\nu(y) \exp\K{(}{)}{ky} &= \exp\K{(}{)}{\psi(k)}\\
&= (1 - k \delta t)^{-\frac{1}{\delta t}}.
\end{align}

Therefore, for $u \in \mathbb{R}_{\geq 0}$
\begin{align}
(\delta u)^{-\frac{1}{\delta t}} &= \int_{\mathbb{R}} \exp\K{(}{)}{\frac{y}{\delta t}} d\nu(y) \exp\K{(}{)}{ - u y}.
\end{align}

From the inverse Laplace transform we obtain:
\begin{align}
d\nu(y) = \frac{(\delta t)^{-\frac{1}{\delta t}}}{\Gamma\K{(}{)}{\frac{1}{\delta t}}}\exp\K{(}{)}{-\frac{y}{\delta t}}y^{\frac{1}{\delta t}}1_{y > 0}\frac{dy}{y}.
\end{align}

Finally
\begin{align}
d\mu(y|x) &= d\nu(y)\exp\K{(}{)}{\phi(x) + \partial_x \phi(x) (y - x)},\\
&= \frac{(\delta t)^{-\frac{1}{\delta t}}}{\Gamma\K{(}{)}{\frac{1}{\delta t}}}\K{(}{)}{\frac{y}{x}}^{\frac{1}{\delta t} }\exp\K{(}{)}{-\frac{1}{\delta t}\K{(}{)}{\frac{y}{x} }}1_{y > 0} \frac{dy}{y}.
\end{align}

\subsubsection{Entropic Non-Compact : Poisson}

Here $\Omega = (0, \infty)$ and
\begin{align}
\sigma : \Omega &\rightarrow \mathbb{R}_{\geq 0},\\
x &\mapsto \sqrt{x}.
\end{align}

Therefore
\begin{align}
\phi(x) = \frac{1}{\delta t} x\log(x) + ax + b.
\end{align}

We choose $x_* = 1$, that is we fix $a,b$ such that $\phi(1) = \partial_x \phi(1) = 0$, hence
\begin{align}
\phi(x) = \frac{1}{\delta t} \K{(}{)}{x\log(x) - x + 1}.
\end{align}

Then
\begin{align}
\psi(k) = \sup_{x \in \Omega} \K{(}{)}{ kx - \phi(x)}.
\end{align}

The maximum is attained where
\begin{align}
k = \partial_x \phi(x) = \frac{1}{\delta t}\log(x).
\end{align}

In particular $\Omega^* = \mathbb{R}$. Solving for $x$,
\begin{align}
x = \exp\K{(}{)}{k \delta t}.
\end{align}

Therefore
\begin{align}
\psi(k) = \frac{\exp\K{(}{)}{k \delta t} - 1}{\delta t}.
\end{align}

Hence
\begin{align}
\int_{\mathbb{R}} d\nu(y) \exp\K{(}{)}{ky} &= \exp\K{(}{)}{\psi(k)}\\
&= \exp\K{(}{)}{\frac{1}{\delta t}\K{(}{)}{e^{k \delta t} - 1}}\\
&= e^{-\frac{1}{\delta t}} \sum_{n \geq 0} \frac{e^{nk\delta t}}{n! (\delta t)^n}.
\end{align}

Therefore
\begin{align}
d\nu(y) &= e^{-\frac{1}{\delta t}} \sum_{n \geq 0} \frac{1}{n! (\delta t)^n} \delta(y - n\delta t)dy.
\end{align}

Finally
\begin{align}
d\mu(y|x) &= d\nu(y)\exp\K{(}{)}{\phi(x) + \partial_x \phi(x) (y - x)}\\
&=   \exp\K{(}{)}{-\frac{x}{\delta t}}\sum_{n \geq 0} \frac{x^n}{n! (\delta t)^n} \delta(y - n\delta t)dy.
\end{align}

Note that the above formula shows that for finite $\delta t$ the process defines a Markov chain on the natural numbers, indeed the state space $\Gamma$ is given by
\begin{align}
\Gamma = \delta t \mathbb{N}_0,
\end{align}

and $0 \in \Gamma$ is the only absorbing state.
\subsubsection{Entropic Compact : Bernoulli}

Here $\Omega = (0, 1)$ and
\begin{align}
\sigma : \Omega &\rightarrow \mathbb{R}_{\geq 0},\\
x &\mapsto \sqrt{x(1 - x)}.
\end{align}

Therefore
\begin{align}
\phi(x) = \frac{1}{\delta t} \K{(}{)}{x \log(x) + (1 - x)\log(1 - x) + ax + b}.
\end{align}

We choose $x_* = q \in (0, 1)$, that is we fix $a,b$ such that $\phi(q) = \partial_x \phi(q) = 0$, hence
\begin{align}
\phi(x) = \frac{1}{\delta t} \K{(}{)}{x\log\K{(}{)}{\frac{x}{q}} + (1 - x)\log\K{(}{)}{\frac{1 - x}{1 - q}}}.
\end{align}

namely $\phi$ is the entropy of the 2-state distribution $(x, 1-x)$ relative to $(q, 1-q)$. The volatility we are considering corresponds to that of a random walk on $(0,1)$ viewed as a space of probability distributions over 2 states.
Then
\begin{align}
\psi(k) = \sup_{x \in \Omega} \K{(}{)}{ kx - \phi(x)}.
\end{align}

The maximum is attained where
\begin{align}
k = \partial_x \phi(x) = \frac{1}{\delta t}\K{(}{)}{\log\K{(}{)}{\frac{x}{q}} - \log\K{(}{)}{\frac{1 - x}{1 - q}} }.
\end{align}

In particular $\Omega^* = \mathbb{R}$. Solving for $x$
\begin{align}
x = \frac{1}{1 + \dfrac{1-q}{q}\exp\K{(}{)}{-k\delta t}}.
\end{align}

Therefore
\begin{align}
\psi(k) = \frac{1}{\delta t}\log\K{(}{)}{(1 - q) + q \exp\K{(}{)}{k \delta t}}.
\end{align}

Hence
\begin{align}
\int_{\mathbb{R}} d\nu(y) \exp\K{(}{)}{ky} &= \exp\K{(}{)}{\psi(k)}\\
&= \K{(}{)}{(1 - q) + q \exp\K{(}{)}{k \delta t}}^{\frac{1}{\delta t}}.
\end{align}

We specialize to the case
\begin{align}
\delta t = \frac{1}{N},
\end{align}

then
\begin{align}
\int_{\mathbb{R}} d\nu(y) \exp\K{(}{)}{ky} &= \exp\K{(}{)}{\psi(k)}\\
&= \sum_{n = 0}^N {N \choose n} q^n (1 - q)^{N - n} \exp\K{(}{)}{\frac{n}{N}k}.
\end{align}

Therefore
\begin{align}
d\mu(y|x) &= \sum_{n = 0}^N {N \choose n} x^n (1 - x)^{N - n} \delta\K{(}{)}{y - \frac{n}{N}}dy.
\end{align}

Note that the above formula shows that for $\delta t = 1/N$ the process defines a Markov chain on a finite set, indeed the state space $\Gamma$ is given by
\begin{align}
\Gamma = \K{\{}{\}}{0, \frac{1}{N}, \dots, \frac{N-1}{N}, 1},
\end{align}
and $0,1 \in \Gamma$ are the only absorbing states.
\subsection{Dimension 2}

Here we construct a 2 dimensional Hessian Martingale generalization of a standard 1 dimensional continuous martingale.
\subsubsection{A Stochastic Volatility Model}

A classic stochastic volatiliy model is given by the SABR model introduced in \cite{SABR}:
\begin{align}
du(t) &= \sigma v(t) u(t)^{\beta} dW_{1, t},\\
dv(t) &= \alpha v(t) dW_{2, t},\\
d\langle W_1, W_2\rangle_t &= \rho dt \s\s \beta \in [0,1], \sigma, \alpha \geq 0, \rho \in [-1,1].
\end{align}

This model, we will show, does not fall into the category of Hessian martingales. On the other hand, we know that there is a Hessian martingale that exactly fits to the marginals of SABR. We will not construct such a martingale, instead we will construct a Hessian martingale model that offers a different stochastic volatility generalization of CEV (Constant Elasticity of Variance) (see \cite{CEV}). The covariance matrix of a Hessian martingale has entries
\begin{align}
C_{uu} &= \frac{\partial_v^2 \phi}{\det H(\phi)},\\
C_{uv} &= -\frac{\partial_v\partial_u \phi}{\det H(\phi)},\\
C_{vv} &= \frac{\partial_u^2 \phi}{\det H(\phi)}.
\end{align}

Those of SABR are given by
\begin{align}
C_{uu}' &= \sigma^2 v^2 u^{2\beta},\\
C_{uu}' &= \sigma \alpha \rho v^2 u^{\beta},\\
C_{vv}' &= \alpha^2 v^2.
\end{align}

If SABR were a Hessian martingale, there would be $\phi(u,v)$ such that
\begin{align}
\partial_u^2 \phi &\propto \frac{1}{v^2},\\
\partial_u\partial_v \phi &\propto \frac{1}{v^2 u^{\beta}}.
\end{align}

Differentiating the first by $v$ and the second by $u$ we see that that would never be the case. On the other hand, we can try to preserve the basic feature that $C_{vv} = \alpha^2 v^2$. Then.
\begin{align}
\partial_v^2 \phi - \frac{(\partial_v\partial_u \phi)^2}{\partial_u^2 \phi} = \frac{1}{\alpha^2 v^2}.
\end{align}

We define
\begin{align}
\phi(u,v) =: \widetilde{\phi}(u,v) - \frac{1}{\alpha^2}\log(v).
\end{align}

It follows
\begin{align}
\det H(\widetilde{\phi}) = 0
\end{align}
This equation is in particular invariant under affine transformations. This will allow us to construct more solutions once we find one. We specialize to the following Ansatz:
\begin{align}
\widetilde{\phi}(u,v) = f(u)g(v).
\end{align}

and consider only the case $g'' \neq 0$\footnote{If $g'' = 0$ we would have that either $f$ or $g$ are constant. The non-trivial case here would be $g$ constant, allowing $f$ to be chosen at will, but that case would just correspond to $u,v$ evolving as two independent 1d processes.}. It follows that
\begin{align}
\frac{f(u)f''(u)}{(f'(u))^2} = \frac{(g'(v))^2} {g(v)g''(v)} = a \in \mathbb{R}.
\end{align}

Let us concentrate on $f$. Then
\begin{align}
\partial_u \log\K{|}{|}{f'(u)} = a \, \partial_u \log\K{|}{|}{f(u)}.
\end{align}

Therefore
\begin{align}
|f'(u)| = e^c |f(u)|^{a}.
\end{align}

Hence:
\begin{align}
f(u) &= \pm(A u + B)^{\frac{1}{1 - a}} \s\s a \neq 1,\\
f(u) &= B e^{Au} \s\s a = 1.
\end{align}

We will concentrate on the case $a \neq 1$. Then, (and using the affine symmetry) the class of $\phi$'s we have found is of the form
\begin{align}
\phi(u,v) = \pm(Au + Bv + C)^{\frac{1}{1 - a}}(Du + Ev + F)^{\frac{a}{a-1}} - \frac{1}{\alpha^2}\log(v) + \mathrm{affine}.
\end{align}

We now define $\beta$ such that
\begin{align}
2 - 2\beta = \frac{1}{1 - a},
\end{align}

and further specialize to the case $A = 1, B = 0, C = 0, E = 0$, $\beta \notin \{1/2, 1\}$\footnote{We can construct the latter as limits, but we will ignore this here.}. Then, requiring convexity, $\phi$ will reduce to the form
\begin{align}
\phi(u,v) = -\frac{1}{(2 - 2\beta)(2\beta - 1)\sigma^2}u^{2 - 2\beta}(v + c)^{2\beta - 1} - \frac{1}{\alpha^2}\log(v) + \mathrm{affine}.
\end{align}

We now compute the Hessian and covariance matrix:
\begin{align}
\partial_u^2 \phi(u,v) &= \frac{1}{\sigma^2}u^{-2\beta}(v + c)^{2\beta - 1},\\
\partial_v\partial_u \phi(u,v) &= -\frac{1}{\sigma^2}u^{1 - 2\beta}(v + c)^{2\beta - 2},\\
\partial_v^2 \phi(u,v) &= \frac{1}{\sigma^2}u^{2 - 2\beta}(v + c)^{2\beta - 3} + \frac{1}{\alpha^2 v^2},\\
\det H(\phi)(u,v) &= \frac{1}{\alpha^2 v^2}\partial_u^2 \phi(u,v).
\end{align}

Therefore
\begin{align}
C_{uu} &= \alpha^2 \K{(}{)}{\frac{u}{v + c}}^2 v^2 + \sigma^2 \K{(}{)}{\frac{u}{v + c}}^{2\beta}(v + c),\\
C_{uv} &= \alpha^2 \K{(}{)}{\frac{u}{v + c}}v^2,\\
C_{vv} &= \alpha^2 v^2.
\end{align}

Notice how, when $\alpha = 0$, the stochastic process reduces to a CEV model as desired. Therefore we can intepret the present model as a generalization of CEV that incorporates stochastic volatility.
The correlation $\rho$ here is not constant. Moreover it is never negative and its square is given by
\begin{align}
\rho^2  &= \frac{1}{1 + \dfrac{\sigma^2}{\alpha^2} \K{(}{)}{\dfrac{u}{v+c}}^{2\beta - 2}(v + c)v^{-2}}.
\end{align}

\end{document}